\documentclass[12pt]{article}

\usepackage[english]{babel}

\usepackage{amsthm}

\usepackage{cite}

\usepackage{enumerate}
\usepackage{amsfonts,amsmath,amsthm,amscd,amssymb,latexsym}

\begin{document}

\newtheorem{theorem}{Theorem}[section]
\newtheorem{lemma}{Lemma}[section]
\newtheorem{proposition}{Proposition}[section]
\newtheorem{corollary}{Corollary}[section]
\newtheorem{definition}{Definition}[section]
\theoremstyle{definition}
\newtheorem{example}{Example}
\newtheorem{remark}{Remark}
\newenvironment{romanlist}{
	\def\theenumi{\roman{enumi}}\def\theenumii{\alph{enumii}}
	\def\labelenumi{(\theenumi)}\def\labelenumii{(\theenumii)}%
	\let\item\Item
	\begin{enumerate}
	}{
	\end{enumerate}}
	
	\let\Item\item
\newenvironment{enumeroman}{%
  \def\theenumi{\roman{enumi}}\def\theenumii{\alph{enumii}}%
  \def\labelenumi{(\theenumi)}\def\labelenumii{(\theenumii)}%
		\let\item\Item
  \begin{enumerate}%
}{%
  \end{enumerate}}
	
\def\address#1{\expandafter\def\expandafter\@aabuffer\expandafter
	{\@aabuffer{\affiliationfont{#1}}\relax\par
	\vspace*{13pt}}}

\setcounter{page}{1}

\markboth{A.~V.~Tushev}{On the Primitive Irreducible Representations of Finitely Generated  
Linear Groups of Finite Rank}



\title{ON THE PRIMITIVE IRREDUCIBLE REPRESENTATIONS OF FINITELY GENERATED 
LINEAR GROUPS OF FINITE RANK}

\author{A.~V.~TUSHEV}

\address{Department of Mathematics, Dnipro National University,\\
Gagarin Avenue 72,Dnipro, 49065, Ukraine\\
tushev@mmf.dnu.edu.ua}

\maketitle

\begin{center}
 
\small{Department of Mathematics, Dnipro National University,\\
Gagarin Avenue 72, Dnipro, 49065, Ukraine\\
tushev@mmf.dnu.edu.ua}

\end{center}

\begin{abstract}
In the paper we study finitely generated linear  groups of finite rank which have faithful irreducible primitive 
representations over a field of characteristic zero. We prove that if an infinite finitely generated 
linear group $G$ of finite rank has a faithful irreducible primitive representation over a field of 
characteristic zero then the $FC$-center  $\Delta (G)$ of $G$ is infinite. 
\end{abstract}




\section{Introduction} 

      A group $G$ is said to have finite (Prufer) rank if there is a positive integer $m$ 
such that any finitely generated subgroup of $G$ may be generated by $m$ elements; 
the smallest $m$ with this property is the rank $r(G)$ of $G$. A group $G$ is said 
to be of finite torsion-free rank if it has a finite series each of whose factor is 
either infinite cyclic or locally finite; the number ${r_0}(G)$ of infinite cyclic 
factors in such a series is the torsion-free rank of $G$. 
\par  
If the group $G$ has a finite series each of whose factor is either cyclic or quasi-cyclic then 
$G$ is said to be minimax; the number $m(G)$ of infinite factors in such a series is the minimax 
length of $G$. If in such a series all infinite factors are cyclic then the group $G$ is said to 
be polycyclic; the number $h(G)$ of infinite factors in such a series is the Hirsch number of $G$. 
It is well known that ${r_0}(G) = m(G) = h(G)$ if the group $G$ is polycyclic.  
\par

	Let $G$ be a group, $R$ be a ring and $RG$ be the group ring. Let $H$ be a subgroup of the group 
$G$ and let $U$ be a right $RH$-module. Since the group ring $RG$ can be considered as a left 
$RH$-module, we can define the tensor product $U{ \otimes _{RH}}RG$, which is a right $RG$-module 
named the $RG$-module induced from the $RH$-module $U$. Moreover, if $M$ is an $RG$-module and $U \le M$, 
then 
\begin{equation}
                      M = U{ \otimes _{RH}}RG                            \label{1.1}
\end{equation}
if and only if 
\begin{equation}
           M = { \oplus _{t \in T}}Ut,                                     \label{1.2}
\end{equation}
where $T$ - is a right transversal to subgroup $H$ in $G$. 
\par   
If a $kG$-module $M$ of some representation  $\varphi $ of a group $G$ over a field $k$ is 
induced from some  $kH$-module $U$, where $H$ is a subgroup of the group $G$, then we say that 
the representation $\varphi $ is induced from a representation $\phi $ of subgroup $H$, 
where $U$ is the module of the representation $\phi $. 
\par  
	Eq. (\ref{1.2}) shows that properties of the $RG$-module $M$ and the $RH$-module 
$U$ are closely related. So, Eq. (\ref{1.1}) may be very useful if properties of 
the $RH$-module $U$ are well studied. For instance, in the case where the group $H$ is 
polycyclic, we have a deeply developed theory (see \cite{Wehr09}). Eq. (\ref{1.1}) 
also may be applied in the case where the group $G$ has finite torsion-free rank ${r_0}(G)$ 
if ${r_0}(H) < {r_0}(G)$ because then we can use the induction on ${r_0}(G)$. However, on 
the place of ${r_0}(G)$ there may be another rank of the group $G$ or the minimax length 
$m(G)$ if $G$ is minimax. 
\par    
	The subgroup $H \le G$ which provides (\ref{1.1}) is said to be a control subgroup for the 
$RG$-module $M$. Various types of control subgroups for modules over group rings of groups 
of finite rank were considered in \cite{Tush95,Tush96,Tush99,Tush2000}. 
\par    
     Let $G$ be a group , let $k$ be a field and let $M$ be a $kG$ -module. The module $M$ is 
said to be primitive if it is not induced from any $kH$-submodule for any subgroup $H < G$. 
Recall that the representation $\varphi $ is said to be faithful if $Ker\varphi  = 1$, it is 
equal to ${C_G}(M) = 1$. 
\par  
	Certainly, primitive irreducible modules are basic subjects for investigations when we are 
dealing with induced modules. Naturally, the following question appears: what can be said 
on the construction of a group $G$ if it has a faithful primitive irreducible representation 
over a field $k$? It should be noted that there are many results which show that the existence 
of a faithful irreducible representation of a group $G$ over a field $k$ may have essential 
influence on the structure of the group $G$ (see for instance 
\cite{Szec16,SzTu17,Tush90,Tush93,Tush2000,Tush02,Tush12}). 
\par   
	In \cite{Harp77} Harper solved a problem raised by Zaleskii and proved that any not 
abelian-by-finite finitely generated nilpotent group has an irreducible primitive representation 
over a not locally finite field. In \cite{Tush02} we proved that if a minimax nilpotent group $G$ 
of class 2 has a faithful irreducible primitive representation over a finitely generated field of 
characteristic zero then the group $G$ is finitely generated. In \cite{Harp80} Harper studied 
polycyclic groups which have faithful irreducible primitive representations. It is well known that any 
polycyclic group is finitely generated soluble of finite rank and meets the maximal condition 
for subgroups (in particular, for normal subgroups). In \cite{Tush2000} we showed that in the 
class of soluble groups of finite rank with the maximal condition for normal subgroups only 
polycyclic groups may have faithful irreducible primitive representations over a field of 
characteristic zero. 
\par   
If $G$ is a group then the $FC$-center $\Delta(G)=\left\{{g\in G|\left|{G:{C_G}(g)}\right|<\infty}\right\}$ 
of $G$ is a characteristic subgroup of $G$. In \cite[Theorem A]{Harp80} Harper proved that if a 
polycyclic group $G$ has a faithful primitive irreducible representation over a field $k$ then 
$\Delta (G)$ is rather large in the sense that $\Delta (G) \cap H > 1$ for any subgroup $1 \ne H$ 
of $G$ such that $|G:{N_G}(H)| < \infty $.  
\par   
By Auslander theorem, any polycyclic group $G$ is linear and it is well known that the group $G$ is 
finitely generated of finite rank. In the presented paper we study finitely generated linear (over a 
field of characteristic zero) groups of finite rank which have faithful irreducible primitive 
representations over a field of characteristic zero. We prove that if an infinite finitely generated 
linear group $G$ of finite rank has a faithful irreducible primitive representation over a field of 
characteristic zero then $\Delta (G)$ is infinite (see Theorem 6.1). 
\par   
Our methods of investigations are based on the following techniques introduced by Brookes in \cite{Broo88} 
for the case of polycyclic groups. Let $N$ be a group and let $K$ be a normal subgroup of $N$ such that 
the quotient group $N/K$ is torsion-free abelian of finite rank. Let $R$ be a ring and let $W$ be a finitely 
generated $RN$-module. Let $I$ be an $N$-invariant ideal of $RK$ such that $\left| {K/{I^\dag }} \right| < 
\infty $ and $k = R/(R \cap I)$ is a field, where ${I^\dag } = G \cap (1 + I)$. Then $N/{I^\dag }$ 
has a central subgroup $A$ of finite index (see Lemma 2.2(i)). So, the quotient module $W/WI$ may be 
considered as a finitely generated $kA$-module. One may ask how deeply properties of the 
$RN$-module $W$ are related with properties of the $kA$-module $W/WI$. However, for our needs the 
approach is quite fruitful and we consider relations between $W$ and a set of prime ideals of $kA$ 
which are minimal over $Ann_{kA} (W/WI)$.


\section{ Auxiliary Results}

Let $A$, $B$ and $H$ be subgroups of a group $G$ such that $B$ is a normal subgroup of $A$. We say that $A$ is an $H$-invariant subgroup if  $H \leqslant {N_G}(A)$ . A section $A/B$ will be called an $H$-invariant section of $G$ if $A$ and $B$ are $H$-invariant subgroups.  
\par  
Let $H$ be a subgroup of a group $G$, the subgroup $H$ is said to be dense in $G$ if for any $g \in G$ there is an integer $n \in \mathbb{N}$ such that ${g^n} \in H$. If ${g^n} \in G\backslash H$ for any $n \in \mathbb{N}$ and any $g \in G\backslash H$ then the subgroup $H$ is said to be isolated in $G$. If the group $G$ is locally nilpotent then the isolator $i{s_G}(H) = \{ g \in G|{g^n} \in H$ for some $n \in \mathbb{N}\} $ of $H$ in $G$ is a subgroup of $G$ and if $H$ is a normal subgroup then so is $i{s_G}(H)$. 
\par  
If $G$ is a group then $G'$ denotes the derived subgroup of $G$. 


\begin{lemma}
Let $D$ be a nilpotent torsion-free normal subgroup of nilpotency class two of a group $G$. Let  $K$ be the centre of $D$ and let $L = i{s_D}(D')$. 
Suppose that there is a $G$-invariant subgroup $P \leqslant D$ such that $L \leqslant P$ , the quotient group $P/L$ is polycyclic, $P \cap K = L$ 
and $PK$ has finite index in $D$. Then: 
\begin{romanlist}
\item $P'$ is a $G$-invariant polycyclic dense subgroup of $L$;
\item for any finitely generated subgroup $H \leqslant G$ the group $D$ has an $H$-invariant normal polycyclic non-abelian subgroup. 
\end{romanlist}
\end{lemma}

\begin{proof} It is well known that in any group $D$ of nilpotency class two the following relation holds for any $a,b,c,d \in D$:
\begin{equation}
                    [ab,cd] = [a,c][a,d][b,c][b,d].                  \label{2.1}     
\end{equation}
\par
	(i) It follows from (\ref{2.1}) that $[{a^m},{b^n}] = {[a,b]^{mn}}$ for any $a,b \in D$ and any $m,n \in \mathbb{N}$. Then, as $PK$ has 
finite index in $D$, it implies that $(PK)'$ is a dense subgroup of $D'$ and hence $(PK)'$ is a dense subgroup of $L = i{s_D}(D')$. Since $K$ 
is the centre of $D$,  it follows from (\ref{2.1}) that $(PK)' = P'$ and hence $P'$ is a dense subgroup of $L$. 
\par  
	As the quotient group $P/L$ is finitely generated, there is a finitely generated  subgroup $A = \left\langle {{a_1},...,{a_t}} \right\rangle  \leqslant P$ such that $P = LA$. Then it easily follows from (\ref{2.1})  that $P' = (LA)' = A' = (\left\langle {{a_1},...,{a_t}} \right\rangle )' = \left\langle {[{a_i},{a_j}]|1 \leqslant i,j \leqslant t} \right\rangle $ and hence $P'$ is a finitely generated subgroup of $P$. Since $P$ is a $G$-invariant subgroup and ${[a,b]^g} = [{a^g},{b^g}]$ for all $a,b \in D$ and $g \in G$, we can conclude that $P'$ is a $G$-invariant subgroup. 
\par  
(ii) At first, we should note that $L$ is a $G$-invariant subgroup as the isolator in $P$ of the $G$-invariant subgroup $P'$. Since $P'$ is a finitely generated dense subgroup of an abelian torsion-free subgroup $L$, we can conclude that $r(L) < \infty $ and hence $L/P'$ is an union of an ascending series of its finite $G$-invariant subgroups. 
\par  
As the quotient group $P/L$ is finitely generated, there is a finitely generated  subgroup $A = \left\langle {{a_1},...,{a_t}} \right\rangle  \leqslant P$ such that $P = LA$. Let $H = \left\langle {{h_1},...,{h_m}} \right\rangle $ be a finitely generated subgroup of $G$. As $P = LA$, for any ${a_i}$ and ${h_j}$ there are ${u_{ij}} \in L$ and ${v_{ij}} \in A$ such that ${a_i}^{{h_j}} = {u_{ij}}{v_{ij}}$, where $1 \leqslant i \leqslant t$ and $1 \leqslant j \leqslant m$. Since $L/P'$ is an union of an ascending series of its finite $G$-invariant subgroups, there is a $G$-invariant subgroup $T \leqslant L$ such that $P' \leqslant T$, $\left| {T/P'} \right| < \infty $ and all ${u_{ij}} \in T$ and therefore, as ${a_i}^{{h_j}} = {u_{ij}}{v_{ij}}$, we see that all ${a_i}^{{h_j}} \in TA$, where $1 \leqslant i \leqslant t$ and $1 \leqslant j \leqslant m$. It implies that $TA$ is an $H$-invariant normal polycyclic non-abelian subgroup of $D$. 
\end{proof}


\begin{lemma}
Let $G$ be a group of finite rank. Then: 
\begin{romanlist}
\item if the group $G$ has a finite subgroup $D$ such that the quotient group $G/D$ is torsion-free abelian then $G$ has a characteristic central torsion-free subgroup $A$ of finite index; 
\item if the group $G$ is finitely generated linear and has no polycyclic subgroups of finite index then $G$ has a finite series $L \leqslant {G_0} \leqslant G$ of normal subgroups such that $\left| G/{G_0} \right| < \infty $, the quotient group ${G_0}/L$ is polycyclic, the subgroup $L$ is torsion-free nilpotent minimax and has no non-abelian torsion-free polycyclic $G$-invariant sections. 
\end{romanlist}
\end{lemma}
\begin{proof} 
	 (i)  As $G'$ is a characteristic subgroup of $G$, it is easy to note that $C = {C_G}(G')$ is a characteristic subgroup of $G$. Then $B = {C^m}$ is a characteristic subgroup of $G$, where $m = \left| D \right|$. Since $\left| {G'} \right| < \infty $ and $C = {C_G}(G')$, we can conclude that $\left| {G:C} \right| < \infty $. As the group $G$ has finite rank, it is easy to note that $\left| {C:{C^m}} \right| < \infty $ and hence we can conclude that $\left| {G:B} \right| < \infty $. For any $a,b \in C$ and any $g \in G$ we have $[ab,g] = {[a,g]^b}[b,g]$ and, as $[a,g] \in G'$ and $b \in C = {C_G}(G')$, we can conclude that 
\begin{equation}
                     [ab,g] = [a,g][b,g].                  \label{2.2}     
\end{equation}
Any element $a \in B = {C^m}$ may be presented in the form $a = {a_1}^m...{a_t}^m$, where ${a_i} \in C$ and $1 \leqslant i \leqslant t$. It follows from (\ref{2.2}) that for any $g \in G$ we have $[a,g] = {[{a_1},g]^m}...{[{a_t},g]^m}$ then, as $m = \left| D \right|$ and $[{a_i},g] \in G' \leqslant D$, where $1 \leqslant i \leqslant t$, we can conclude that $[a,g] = 1$. It implies that $B$ is a characteristic central subgroup of finite index in $G$. Since the quotient group $G/D$ is torsion-free, we see that so is $B/(B \cap D)$ and it easily implies that $A = {B^m}$ is a torsion-free characteristic central subgroup of finite index in $G$.
\par  
(ii) Since the group $G$ is linear over a field of characteristic zero and $r(G) < \infty $, it follows from Titts theorem \cite[Theorem 10.17]{Wehr73} that the group $G$ is soluble-by-finite. Then it follows from Kolchin-Maltsev theorem \cite[Theorem 3.6]{Wehr73} that the group $G$ has a finite series of normal subgroups ${G_1} \leqslant {G_0} \leqslant G$ such that $\left| G/{G_0} \right| < \infty $, the quotient group ${G_0}/{G_1}$ is abelian and the subgroup ${G_1}$ is torsion-free nilpotent. As the group $G$ is finitely generated, by \cite[Proposition 1.6.11]{Robi96}, so is ${G_0}$ and hence the quotient group $G/{G_1}$ is polycyclic-by-finite. By \cite[Proposition 5.2.8]{LeRo04}, the group ${G_0}$ is minimax. As $r(G) < \infty $, ${G_1}$ has a $G$-invariant subgroup $L$ such that the quotient group ${G_1}/L$ is infinite polycyclic of maximal possible torsion-free rank. It easily implies that $L$ has no $G$-invariant subgroup $Y$  such that the quotient group $N/Y$ is infinite polycyclic and, evidently, the quotient group ${G_0}/L$ is polycyclic. 
\par  
Suppose that $L$ has a non-abelian torsion-free polycyclic $G$-invariant sections $A/B$. The section $A/B$ is acted by $G$ via conjugations and hence $L/{C_L}(A/B) \leqslant Aut(A/B)$, where ${C_L}(A/B)$ is a $G$-invariant subgroup of $L$. Then it follows from \cite[Theorem 21.3.2,Theorem 21.2.1]{KaMe79} that the quotient group $L/{C_L}(A/B)$ is polycyclic. Evidently, ${C_A}(A/B) = {C_L}(A/B) \cap A = C$, where $C/B$ is the centre of $A/B$. As the centre $C/B$ of the torsion-free nilpotent non-abelian group $A/B$ is isolated, we can conclude that $\left| {A/C} \right|  = \infty $ and it easily implies that $\left| {L/{C_L}(A/B)} \right| = \infty $. But it is impossible, because $L$ has no $G$-invariant subgroup $Y$ such that the quotient group $L/Y$ is infinite polycyclic. Thus, we can conclude that $L$ has no non-abelian torsion-free polycyclic $G$-invariant sections. 
\end{proof}

	Let $A$ be an abelian additive torsion-free group acted by a group $G$. Let $\bar A = A \otimes_{\mathbb {Z}} \mathbb{Q}$, we denote by $Soc_G \bar A$ the socle of the $\mathbb{Q}G$-module $\bar A$ and put $Soc_G A = Soc_G \bar A \cap A$. It is not difficult to show that $Soc_G A$ is an isolated $G$-invariant subgroup of $A$. We denote by ${\Delta _G}(A)$ the set of elements of $A$ which have finite orbits under the action of $G$, then ${\Delta _G}(A)$ is a $G$-invariant isolated subgroup of $A$.


\begin{lemma}
Let  $A$ be an abelian torsion-free group acted by a group $G$ and let $H$ be a normal subgroup of finite index in $G$. Then: 

\begin{romanlist}
\item $Soc_H A$ is an isolated $G$-invariant subgroup of $A$ such that ${\Delta _G}(A) \leqslant Soc_H A$ and  ${\Delta_G}(Soc_H A/{\Delta_G} (A))$ is trivial; 
\item ${C_S}(Soc_H A)$ is a normal subgroup  of  $G$ for any normal subgroup $S$ of $G$;  
\item the quotient group $G/{C_G}(Soc_H A)$ is abelian-by-finite if the group $A$ has finite rank and the group $G$ is soluble-by-finite. 
\end{romanlist}
\end{lemma}

\begin{proof} It easy to note that the assertions hold for $A$ if and only if they hold for $\bar A = A \otimes _{\mathbb{Z}}\mathbb{Q}$. So, we can assume that $\bar A = A$ and we consider $A$ as a $\mathbb{Q}G$-module. 
\par  
(i) We show that $Soc_H A$ is a $\mathbb{Q}G$-submodule of $A$. $Soc_H A = { \oplus _{i \in I}}{A_i}$, where ${A_i}$ are simple $\mathbb{Q}H$-modules. As $H$ is a normal subgroup of $G$, ${B_i}g\mathbb{Q}H = {B_i}\mathbb{Q}(gH{g^{ - 1}})g = {B_i}\mathbb{Q}Hg = {B_i}g$ for any $\mathbb{Q}H$-submodule ${B_i}$ of ${A_i}$, any $g \in G$ and any $i \in I$. It implies that ${A_i}g$ is a simple $\mathbb{Q}H$-module and hence ${A_i}g \leqslant Soc_H A$ for any $g \in G$ and any $i \in I$. Therefore, $(Soc_H A)g \leqslant Soc_H A$ for any $g \in G$ and hence $Soc_H A$ is a $G$-invariant subgroup of $G$. 
\par  
It follows from the definition of ${\Delta_G}(A)$ that the quotient group $H/{C_H}(a\mathbb{Q}H)$ is finite for any $a \in {\Delta_G}(A)$. Then, by Maschke’s Theorem, $a\mathbb{Q}H$  is a semi-simple $\mathbb{Q}H$-module for any $a \in {\Delta_G}(A)$. Therefore, ${\Delta_G}(A) \leqslant Soc_H A$. Since $Soc_H A$ is a semi-simple $\mathbb{Q}H$-module, we have $Soc_H A = {\Delta_G}(A) \oplus V$, where $V$ is a $\mathbb{Q}H$-module such that  $V \cong Soc_H A/{\Delta_G}(A)$. As $|G/H|<\infty$, it is easy to note that $\Delta_{G}(A)= 
\Delta_{H}(A)$, and hence $\Delta_{H}(V)=0$. Therefore, $\Delta_{H}(Soc_H A/{\Delta_G} (A))={\Delta_G}(Soc_H A/{\Delta_G} (A))=0$. 
 \par  

	(ii) By (i), $Soc_H A$ is a $G$-invariant subgroup of $G$ and hence for any $h \in {C_S}(Soc_H A)$, any $g \in G$ and any $a \in Soc_H A$ we have $a{h^g} = (a{g^{ - 1}})hg = (a{g^{ - 1}})g = a$ because $a{g^{ - 1}} \in Soc_H A$. Thus, for any $h \in {C_S}(Soc_H A)$ and any $g \in G$ we have ${h^g} \in {C_S}(Soc_H A)$ and the assertion follows. 
	\par
	
(iii) It follows from the definition of $G$ that $H$ has a soluble $G$-invariant subgroup $S$ of finite index in $G$. As the group $A$ has finite rank, ${\dim _\mathbb{Q}}(Soc_H A) < \infty $ and hence by Clifford Theorem \cite[Theorem 1.7]{Wehr73} $S/{C_S}(Soc_H A)$ is a completely reducible soluble subgroup of $GL(n,\mathbb{Q})$. By \cite[Theorem 3.5]{Wehr73}, the quotient group $S/{C_S}(Soc_H A)$ is abelian-by-finite.  By (ii), ${C_S}(Soc_H A)$ is a $G$-invariant subgroup and, as $|G/S|<\infty$, we see that the quotient group $G/C_{S}(Soc_H A)$ is abelian-by-finite.  As ${C_S}(Soc_H A) \leqslant {C_G}(Soc_H A)$, the assertion follows.  
\end{proof}

Let $A$ be an abelian torsion-free group of finite rank acted by a group $G$. The group $G$ acts on $A$ rationally irreducible if $\bar A = A{ \otimes _\mathbb{Z}}\mathbb{Q}$ is an irreducible $\mathbb{Q}G$-module and $A$ is said to be a $G$-plinth if any subgroup of finite index in $G$ acts on $A$ rationally irreducible. If $A$ has a finitely generated dense subgroup $B$ such that $B =  \oplus _{i = 1}^n{B_i}$, where $i{s_A}({B_i})$ is a $G$-plinth for each $1 \leqslant i \leqslant n$, then we say that $A$ is a $G$-polyplinth. A $G$-polyplinth is said to be Noetherian if it satisfies the maximal condition for $G$-invariant subgroups. 


\begin{lemma}
Let $A$ be an abelian torsion-free group of finite rank acted by a soluble-by-finite group $G$. Then:
\begin{romanlist}
\item $Soc_H A = Soc_G A$ for any normal subgroup $H$ of finite index in $G$;
\item $G$ has a normal subgroup H of finite index such that $Soc_H A = Soc_G A$ contains a dense $G$-invariant Noetherian $H$-polyplinth;
\item if $Soc_G A = A$ and $A$ satisfies the maximal condition for $G$-invariant subgroups then ${\Delta _G}(A)$ is a finitely generated isolated $G$-invariant subgroup of $A$; 
\item if $Soc_G A = A$ and $A$ has an isolated $G$-invariant subgroup $B$ such that the quotient group $A/B$ is polycyclic then $A$ has a $G$-invariant polycyclic subgroup $P$ such that $P \cap B = 0$ and $P + B$ has finite index in $A$. 
\end{romanlist}
\end{lemma}

\begin{proof} It easy to note that the assertions (i) and (ii) hold for $A$ if and only if they hold for $\bar A = A{ \otimes _\mathbb{Z}}\mathbb{Q}$. So, in the proof of assertions (i) and (ii) we can assume that $\bar A = A$ and we consider $A$ as a $\mathbb{Q}G$-module. 
\par  
(i) As $dim_\mathbb{Q} A$ is finite, it follows from Clifford Theorem \cite[Theorem 1.7]{Wehr73} that $Soc_X A \geqslant Soc_G A$ for any normal subgroup $X$ of $G$. Since  $dim_\mathbb{Q} A$ is finite,  it implies that there is a $G$-invariant subgroup $H$ of finite index in $G$ such that $Soc_H A = Soc_S A$ for any $G$-invariant subgroup $S$ of finite index in $H$. Then it is sufficient to show that $Soc_H A = Soc_G A$. By Lemma 2.3(i,iii),  $Soc_H A$ is a $G$-invariant subgroup of $A$ and, the quotient group $G/{C_G}(Soc_H A)$ is abelian-by-finite. Therefore, changing $G$ by $G/C_G (Soc_H A)$ we may assume that $G$ is abelian–by-finite. So, $G$ has a normal abelian subgroup $X$ of finite index and changing $H$ by $H \cap X$ we may assume that $H$ is abelian.
\par  
	It is sufficient to show that $Soc_H A$ is a semi-simple $\mathbb{Q}G$-module. Since  $H$ is abelian, we see that $Soc_H A =  \oplus _{i = 1}^k B_i$, where each $B_i$ is a simple $\mathbb{Q}H$-module which is annihilated by a maximal  ideal $P_i$ of $\mathbb{Q}H$. Since $Soc_H A = \sum\nolimits_{i = 1}^k {{B_i}\mathbb{Q}G} $, it is sufficient to show that each ${B_i}\mathbb{Q}G$ is a semi-simple $\mathbb{Q}G$-module. To simplify the denotations, we fix some index $i$ putting ${B_i} = B$ and ${P_i} = P$. Let $K = {N_G}(P)$, evidently $H \leqslant K \leqslant G$ and hence $K$ has finite index in $G$. Since $P$ annihilates $B\mathbb{Q}K$, we can consider $B\mathbb{Q}K$ as a $\mathbb{Q}K/P\mathbb{Q}K$-module. It is easy to note that the quotient ring $\mathbb{Q}K/P\mathbb{Q}K$ is a cross product $F * T$ of the field $F = \mathbb{Q}H/P$ and the group $T = K/H$. So, $B\mathbb{Q}K$ is an $F * T$-module. It follows from \cite[Theorem 4.1]{Pass89} that the cross product $F*T$ is a semi-simple ring and hence so is $\mathbb{Q}K/P\mathbb{Q}K$. As $P$ annihilates $B\mathbb{Q}K$, it implies that  $B\mathbb{Q}K$ is a semi-simple $\mathbb{Q}K$-module. Thus, $B\mathbb{Q}K = \oplus _{i = 1}^r{C_i}$, where each ${C_i}$ is a simple $\mathbb{Q}K$-module which is annihilated by $P$.
\par  
	We have to show that $B\mathbb{Q}G$ is a semi-simple $\mathbb{Q}G$-module. Since $B\mathbb{Q}K =  \oplus _{i = 1}^r{C_i}$, we have $B\mathbb{Q}G = \sum\nolimits_{i = 1}^r {{C_i}} \mathbb{Q}G$ and hence it is sufficient to show that ${C_i}\mathbb{Q}G$ is a simple $\mathbb{Q}G$-module for each $i$. 
To simplify the denotations we fix some index $i$ and put $C=C_i$. 

	Let $T = \{ {t_1},...{t_l}\} $ be a  transversal to $K$ in $G$ then $C\mathbb{Q}G  =  \sum\nolimits_{i = 0}^l {(C\mathbb{Q}K){t_i}} $. As $K = {N_G}(P)$, 
each $\mathbb{Q}H$-submodule $(C\mathbb{Q}K){t_i}$ is annihilated by $P_i  = P^{t_i}$, where $P_i  = P^{t_i}$ are maximal ideals of $\mathbb{Q}H$ such that $P_i  \ne {P_j}$ if $i \ne j$. 
\par   
	So, we have to show that ${C}\mathbb{Q}G$ is a simple $\mathbb{Q}G$-module.  As $C$ is a simple $\mathbb{Q}K$-module, it is sufficient to show that  $C \cap W \neq 0$  for any $\mathbb{Q}G$-submodule $W \ne 0$ of $C\mathbb{Q}G$. Let $0 \ne a \in W$ then $0 \ne a = \sum\nolimits_{i = 1}^l {{a_i}} $, where ${a_i} \in (C\mathbb{Q}K){t_i}$ and ${a_i}$ is annihilated by ${P_i}$. If ${a_i} \ne 0$ and ${a_j} \ne 0$, where $i \ne j$, then, as ${P_i} \ne {P_j}$, there is an element $x \in {P_j}\backslash {P_i}$ and hence  ${a_i}x \ne 0$ but ${a_j}x = 0$. Then the sum  $ax = \sum\nolimits_{i = 1}^l {{a_i}} x$, where ${a_i}x \in (C\mathbb{Q}K){t_i}$, has less non-zero summands than $0 \ne a = \sum\nolimits_{i = 1}^l {{a_i}} $. So, in such way we can decries the number of non–zero summands from $(C\mathbb{Q}K){t_i}$ in $0 \ne a = \sum\nolimits_{i = 1}^l {{a_i}} $ until we obtain $0 \ne b \in (C\mathbb{Q}K){t_i} \cap W$ for some $i$. Then $0 \ne b{t_i}^{ - 1} \in (C\mathbb{Q}K) \cap W$. 
\par
	(ii) By (i), $Soc_H A = Soc_G A = S$ for any normal subgroup $H$ of finite index in $G$, we denote by $n(H)$ the number of irreducible direct summands in $Soc_H A = S$. Since the numbers $n(H)$ are bounded by $r(A)$, there is a normal subgroup $H$ of finite index in $G$ with maximal possible $n(H)$. Let $Soc_H A = S =  \oplus _{i = 1}^n{S_i}$ , where ${S_i}$ are irreducible $\mathbb{Q}H$-modules, $1 \leqslant i \leqslant n$ and $n = n(H)$. Suppose that there is a submodule ${S_i}$ which is not $\mathbb{Q}L$-irreducible for some subgroup $L$ of finite index in $H$. It is well known that there is a normal subgroup $T$ of finite index in $G$ such that $T \leqslant L$ and hence ${S_i}$ is not $\mathbb{Q}T$-irreducible. It easily implies that ${S_i}$ has a proper $\mathbb{Q}T$-irreducible submodule and hence $n(T) > n(H)$ but it contradicts to the choice of the subgroup $H$. Thus, each ${S_i}$ is an irreducible $\mathbb{Q}L$-module for any subgroup $L$ of finite index in $H$ and hence $Soc_H A = S$ is an $H$-polyplinth. 
\par  
	Let $Soc_H A = Soc_G A = S =  \oplus _{i = 1}^k{A_i}$, where ${A_i}$ is  an irreducible $\mathbb{Q}G$-module, $1 \leqslant i \leqslant k$. Let $0 \ne {b_i} \in {A_i}$, ${B_i} = {b_i}\mathbb{Z}G$ and $B =  \oplus _{i = 1}^k{B_i}$. As ${A_i}$ is an irreducible $\mathbb{Q}G$-module, ${A_i} = {b_i}\mathbb{Q}G$ and hence $B =  \oplus _{i = 1}^k{B_i}$ is a dense subgroup of $Soc_G A = S =  \oplus _{i = 1}^k{A_i}$. Besides, as ${A_i}$ is an irreducible $\mathbb{Q}G$-module, we see that for any $0 \ne x \in {B_i}$ there is $m \in \mathbb{N}$ such that ${b_i}m \in x\mathbb{Z}G$ and hence ${B_i}m \leqslant x\mathbb{Z}G$. Since $\left| {{B_i}/{B_i}m} \right| < \infty $, we can conclude that $\left| {{B_i}/x\mathbb{Z}G} \right| < \infty $ for any $0 \ne x \in {B_i}$ and it easily implies that ${B_i}$ is a Noetherian $\mathbb{Z}G$-module for each $1 \leqslant i \leqslant n$ and hence so is $B =  \oplus _{i = 1}^k{B_i}$. It follows from \cite[Theorem A]{Wils70} that$B =  \oplus _{i = 1}^k{B_i}$ meets the maximal condition for $H$-invariant subgroups. As $Soc_H A = S$ is an $H$-polyplinth and $B =  \oplus _{i = 1}^k{B_i}$ is dense in $Soc_H A = S$, we can conclude that $B =  \oplus _{i = 1}^k{B_i}$ is a $G$-invariant dense in $S$  Noetherian $H$-polyplinth. 
\par  
(iii) Since $r(A) < \infty $, it follows from the definition of ${\Delta _G}(A)$ that it has a finitely generated dense $G$-invariant subgroup $D$. As $A$ satisfies the maximal condition for $G$-invariant subgroups, it implies that $\left| {{\Delta_G}(A)/D} \right| < \infty $ and hence ${\Delta_G}(A)$ is a finitely generated subgroup.  
\par  
(iv) Let $\bar A = A{ \otimes _\mathbb{Z}}\mathbb{Q}$ and $\bar B = B{ \otimes _\mathbb{Z}}\mathbb{Q}$. Since $Soc_G A = A$, we see that $\bar A$ is a semi-simple $\mathbb{Q}G$-module and hence there is a submodule $\bar P$ of $\bar A$ such that $\bar P \cap \bar B = 0$ and $\bar P + \bar B = \bar A$. Put $P = \bar P \cap A$ then $P \cap B = 0$ and $P + B$ is dense in $A$. As the quotient group $A/B$ is polycyclic, we can conclude that $P$ is a $G$-invariant polycyclic subgroup of $A$ and $P + B$ has finite index in $A$.  
\end{proof}

Let $R$ be a ring, an $R$-module is said to be uniform if any its two non-zero submodules have non-zero intersection. 


\begin{lemma}
Let $G$ be a group, let $N$ be a normal subgroup of $G$. Let $R$ be a domain and let $0 \ne W$ be an $RG$-module. Suppose that the subgroup $N$ is torsion-free nilpotent. If $An{n_{SN}}w \ne 0$ for some $0 \ne w \in W$ and some subring $S \leqslant R$ then the submodule $wRG$ is $SN$-torsion. 
\end{lemma}

\begin{proof} Any element $0 \ne x \in wRG$ may be presented in the form $x = {\sum\nolimits_{i = 1}^m {w{d_i}g} _i}$, where ${d_i} \in R$, ${g_i} \in G$ and $m \in \mathbb{N}$. Put $J = An{n_{SX}}w \ne 0$, it follows from \cite[Corollary 37.11]{Pass89} that the ring $SN$ is an Ore domain and hence $SN$ is uniform. Therefore, $I =  \cap _{i = 1}^mg_i^{ - 1}J{g_i} \ne 0$. Evidently, $I \leqslant An{n_{SN}}(x)$ and hence the submodule $wRG$ is $SN$-torsion.
\end{proof}


\section{Large Culling Ideals of Group Rings of Minimax Nilpotent Groups}

Let $K$ be a normal subgroup of a group $L$ such that the quotient group $L/K$ is free abelian with free generators $K{x_1},K{x_2}, ...,K{x_n}$. We say that $\chi  = \{ \left\langle {K,\{ {x_j}|j \in J\} } \right\rangle |J \subseteq \left\{ {1,..., n} \right\}\} $ is a full system of subgroups of $L$ over $K$.
\par  
Let $k$ be a field, let $V$ be a $kK$-module and let $W$ be an image of $V{ \otimes _{kK}}kL$ under a $kL$-module homomorphism $\alpha $. Put $$\chi (W)  =  \left\{ {X \in \chi |\ker \alpha  \cap (V{ \otimes _{kK}}kX} \right\} = 0\} $$ 
and $M\chi (W) $ denotes the set of maximal elements of $\chi (W) $. 


\begin{lemma}
 Let $L$ be a minimax nilpotent group and $K$ be a normal subgroup of $L$ such that the quotient group $L/K$is free abelian with free generators $K{x_1}, K{x_2}, ..., K{x_n}$. Let $\chi  = \{ \left\langle {K,\{ {x_j}|j \in J\} } \right\rangle |J \subseteq \left\{ {1,...,n} \right\}\} $ be a full system of subgroups of $L$ over $K$. Then the group $L$ has a finitely generated dense subgroup $H$ such that if a subgroup ${H_1} \leqslant L$ contains $H$ then  for any $L$-invariant subgroup ${I^\dag } \leqslant K$ of  finite index in $K$ 
\begin{romanlist}
\item $K = {I^\dag }{K_1}$  and  $L = {I^\dag }{H_1}$, where ${K_1} = K \cap {H_1}$;
\item $\chi  =\left\{   {I^\dag} X|X \in {\chi_H}  \right\}$, 
where ${\chi _H} = \{ \left\langle {{K_1},\{ {x_j}|j \in J\} } \right\rangle |J \subseteq \left\{ {1,...,n} \right\}\} $ 
is a full system of subgroups of ${H_1}$ over ${K_1} = K \cap {H_1}$.
\end{romanlist}
\end{lemma}

	\begin{proof} Since the group $L$ is nilpotent and minimax, the subgroup $K$ has a finite central series 
\begin{equation}
 \left\langle 1 \right\rangle  = {C_0} \leqslant {C_1} \leqslant ... \leqslant {C_n} = K                                \label{3.1}
\end{equation}
all of whose quotients are cyclic or quasi-cyclic. Let $\left\{ {{y_1},...,{y_m}} \right\}$ be a set of elements of $K$  which generate all cyclic quotients of (\ref{3.1}) and let $H$ be a subgroup of $L$ generated by $\left\{ {{x_1},{x_2},...,{x_n}} \right\} \cup \left\{ {{y_1},...,{y_m}} \right\}$. Then $H$ is a finitely generated dense subgroup of $L$. Let ${H_1}$ be a subgroup of $L$ such that $H \leqslant {H_1}$. 
\par  
(i) Since $\{ {y_1},...,{y_m}\}  \subseteq {K_1} = K \cap {H_1}$, all non-trivial quotients of the series 
             $$ {K_1} = {C_0}{K_1} \leqslant {C_1}{K_1} \leqslant ... \leqslant {C_n} = K $$                                      
are quasi-cyclic and hence divisible.  Let ${I^\dag }$ be an $L$-invariant subgroup  of  finite index in $K$, as ${K_1}$ meets the conditions of  \cite[Lemma 2.1.7]{Tush2000}, we have $K = {I^\dag }{K_1}$. As $\{ {x_1},{x_2},...,{x_n}\}  \subseteq {H_1}$, we have $L = K{H_1}$ and hence $L = K{H_1} = ({I^\dag }{K_1}){H_1} = {I^\dag }({K_1}{H_1}) = {I^\dag }{H_1}$.  
\par  
	(ii) It follows from the definitions of $\chi $ and ${\chi _H}$ that $\chi  =\left\{ {KX|X \in {\chi _H}} \right\}$. Then, as $K = {I^\dag }{K_1}$ and ${K_1} \leqslant X$ for any $X \in {\chi _H}$, we have $\chi  = \left\{ {KX|X \in {\chi _H}} \right\} = \left\{ {({I^\dag }{K_1})X|X \in {\chi _H}} \right\} = \left\{ {{I^\dag }({K_1}X)|X \in {\chi _H}} \right\} = \left\{ {{I^\dag }X|X \in {\chi _H}} \right\}$.
\end{proof}


\begin{lemma}
 Let $N$ be a minimax nilpotent torsion-free group and let $K$ be a normal subgroup of  $N$ such that the quotient group $N/K$is torsion-free abelian. Let $L$ be a normal subgroup of $N$ such that $K \leqslant L$ and the quotient group $L/K$ is free abelian. Let $\chi $ be a full system of subgroups of $L$ over $K$. Let $k$ be a field and let $0 \ne W$ be a $kN$-module which is $kK$-torsion-free. Then there is a $kN$-submodule $0 \ne V \leqslant W$ such that for any element $0 \ne a \in V$:
\begin{romanlist}
\item there is a finitely generated dense subgroup $H \leqslant L$ such that $akL  =  akH{ \otimes _{kH}}kL$;
\item the module $V$ is $kX$-torsion-free for any $X \in \chi (akL)$ and $kX$-torsion for any $X \in \chi \backslash \chi (akL)$;
\item $\chi (akL) = \chi (bkL)$ for any element $0 \ne b \in V$;
\item in the case where $N = L$, for any subgroup $X \in M\chi (akL)$ we have $bkL \cap ak(X \cap H) \ne 0$ for any element $0 \ne b \in akL$ and any finitely generated dense subgroup $H \leqslant L$. In particular, $bkL \cap akH \ne 0$, $bkL \cap akX \ne 0$ and the module $akL$ is uniform. 
\end{romanlist}
\end{lemma}

\begin{proof} It is easy to note that in the proof the module $W$ may be changed by any its proper submodule. 
\par  
(i). By \cite[Lemma 2.2.5]{Tush2000}, the assertion holds for any element $0 \neq a \in W $. 
\par  
 (ii). Let $X \in \chi $ if the module $W$ is $kX$-torsion-free then $X \in \chi (akL)$ for any element $0 \ne a \in W$. If $An{n_{kX}}w \ne 0$ for some $0 \ne w \in W$ then we may change $W$ by $wkN$. As $X$ is a normal subgroup of $N$, it follows from Lemma 2.5 that the module $V = wkN$ is $kX$-torsion and hence $X \in \chi \backslash \chi (akL)$ for any element $0 \ne a \in V$. Taking $Y \in \chi \backslash \{ X\} $ and repeating the above arguments we obtain a $kN$-submodule $0 \ne {V_1} \leqslant V$ which is either $kY$-torsion or $kY$-torsion-free. Continuing this process we see that it is terminated because the set $\chi $ is finite.
\par  
 (iii). The assertion follows from (ii).
\par  
(iv). By (ii), we can assume that for any $0 \ne w \in W$ and $X \in M\chi (wRL)$ the $kN$-module $wkN$ is $kX$-torsion-free and hence, by \cite[Lemma 2.2.4(ii)]{Tush2000}, there exists a $kL{(kX)^{ - 1}}$-module $wkL{(kX)^{ - 1}}$. Then it easily follows from the maximality of $X$ that $wkL{(kX)^{ - 1}}$ has finite dimension over the field $kX{(kX)^{ - 1}}$. Therefore, we can choose the element $0 \ne v \in wkL$ such that $vkL{(kX)^{ - 1}}$ is a simple $kL{(kX)^{ - 1}}$-module. Put ${V_1} = vkL$. Then for any $0 \ne a,b \in V$ we have $V = akL{(RX)^{ - 1}} = bkL{(RX)^{ - 1}}$ and hence $a \in bkL{(kX)^{ - 1}}$, it easily implies that $bkL \cap akX \ne 0$. Since $H \cap X$ is a dense subgroup of $X$, it follows from \cite[Lemma 2.2.3(ii)]{Tush2000} that $bkL \cap ak(H \cap X) \ne 0$. The last equations also shows that $bkL \cap akH \ne 0$ and the module ${V_1} = vkL$ (and hence any its proper submodule) is uniform. Thus, we obtained a submodule $0 \ne {V_1} = vkL$ such that the assertion (iv) holds for the chosen subgroup $X \in M\chi (vRL)$ and elements of ${V_1}$. 
\par  
Suppose that there is a subgroup $Y \in M\chi (vRL)\backslash \{ X\} $. Applying the above arguments to $Y$ and the submodule ${V_1}$ we may obtain a submodule $0 \ne {V_2} \leqslant {V_1}$ such that the assertion (iv) holds for the chosen subgroups $X,Y \in M\chi (vRL)$ and elements of ${V_2}$. Continuing this process we see that it is terminated because the set $\chi $ is finite. 
\end{proof}

Let $G$ be a group and let $K$ be a normal subgroup of $G$. Let $R$  be a ring and let $I$ be a $G$-invariant ideal of the group ring $RK$ then ${I^\dag } = G\bigcap {(I + 1)} $ is a $G$-invariant subgroup of $G$. We say that the ideal $I$ is $G$-large if ${R \mathord{\left/
 {\vphantom {R {(R\bigcap I )}}} \right.
 \kern-\nulldelimiterspace} {(R\bigcap I )}} = k$ is a field, $\left| {{K \mathord{\left/
 {\vphantom {K {{I^\dag }}}} \right.
 \kern-\nulldelimiterspace} {{I^\dag }}}} \right| < \infty $ and $I = (RF\bigcap I )RK$, where $F$ is a $G$-invariant subgroup of $K$ such that ${I^\dag } \leqslant F$ and the quotient group ${F \mathord{\left/
 {\vphantom {F {{I^\dag }}}} \right.
 \kern-\nulldelimiterspace} {{I^\dag }}}$ is abelian. If $R$ is a field then, certainly, $R = k$. If $N$ is a subgroup of $G$ such that $K \leqslant N$ and ${F \mathord{\left/
 {\vphantom {F {{I^\dag }}}} \right.
 \kern-\nulldelimiterspace} {{I^\dag }}}$ is a central section of $N$ then we say that the ideal $I$ is $N$-central. 


\begin{lemma}
Let $K$ be a normal subgroup of a torsion-free minimax nilpotent  group $N$ such that the quotient group $N/K$ is torsion-free abelian. Let $L$ be a dense subgroup of $N$ such that $K \leqslant L$ and  the quotient group $L/K$ is free abelian. Let $k$ be a field, let $I$ be an $N$-large ideal of $kK$ and let $A \leqslant L/{I^\dag }$ be a central dense free  subgroup of  $N/{I^\dag }$. Let $0 \ne W$ be a $kN$-module which is $kK$-torsion-free. Then there is a $kL$-submodule $0 \ne V \leqslant W$ such that for any element $0 \ne a \in V$ and any element $0 \ne b \in akN$: 
\begin{romanlist}

	\item $bkL$ is not isomorphic to any proper section of $agkL$ for any $g \in N$; 
	\item the $kA$-module $bkL/bkLI$ has a finite series each of whose quotient is isomorphic to some section of the $kA$-module $akL/akLI$; 
\item there is a $kL$-submodule $0 \ne ckL \leqslant akL$ such that the $kA$-module $ckL/ckLI$ has a finite series each of whose quotient is isomorphic to some section of the $kA$-module $bkL/bkLI$.
\end{romanlist}
\end{lemma}

	\begin{proof} (i) By Lemma 3.2(ii), we can assume that for any element $0 \ne a \in W$ and any subgroup $X \in M\chi (akL)$ the module $W$ is $kX$-torsion-free. By Lemma 3.2(iv), there is a $kL$-submodule $0 \ne V \leqslant W$ such that for any element $0 \ne a \in V$ and any element $0 \ne b \in akL$ we have $bkL \cap akX \ne 0$. Therefore, for any submodule $0 \ne U \leqslant akL$ the quotient module $akL/U$ is $kX$-torsion. As $K \leqslant X$ and the quotient group $N/K$ is abelian, $X$ is a normal subgroup of $N$ and hence $agkL/U$ is $kX$-torsion for any $kL$-submodule $0 \ne U \leqslant agkL$ and any element $g \in N$. It easily implies that any proper section of $agkL$ is $kX$-torsion. On the other hand, the module $akN$ is $kX$-torsion-free and the assertion follows. 
\par  
(ii). Step 1. At first, we prove the assertion in the case, where $0 \ne b \in akL$. By \cite[Lemma 2.2.5]{Tush2000}, there is a finitely generated dense subgroup $H \leqslant L$ such that $akL = akH{ \otimes _{kH}}kL$. Evidently, the subgroup $H$ can be chosen such that $0 \ne b \in akH$ and ${I^\dag }H = L$ then $bkL = bkH{ \otimes _{kH}}kL$. 
\par  
Put ${K_1} = K \cap H$ and ${I_1} = I \cap k{K_1}$ then it follows from \cite[Lemma 2.2.6]{Tush2000} that $akL/akLI \cong akH/akH{I_1}$ and $bkL/bkLI \cong bkH/bkH{I_1}$ considered as $kA$-modules, where $A$ is a central subgroup of finite index in $L/{I^\dag }$. Therefore, it is sufficient to show that the $kA$-module $bkH/bkH{I_1}$ has a finite series each of whose quotient is isomorphic to some section of 
the $kA$-module $akH/akH{I_1}$.  
\par  
Since the ideal ${I_1}$ is $H$-large, there exists an $H$ -invariant subgroup  $F \leqslant {K_1}$ such that the quotient group ${F \mathord{\left/
 {\vphantom {F {I_1^\dag }}} \right.
 \kern-\nulldelimiterspace} {I_1^\dag }}$ is abelian and ${I_1} = (kF{\bigcap I _1})k{K_1}$. As $\left| {{F \mathord{\left/
 {\vphantom {F {I_1^\dag }}} \right.
 \kern-\nulldelimiterspace} {I_1^\dag }}} \right| < \infty $, we have $\left| {H:C} \right| < \infty $ , where $C = {C_H}({F \mathord{\left/
 {\vphantom {F {I_1^\dag }}} \right.
 \kern-\nulldelimiterspace} {I_1^\dag }})$. The arguments of the proof of  \cite[Lemma 2.4.1(i)]{Tush2000} shows that $(kF{\bigcap I _1})kC$ is a polycentral ideal of $kC$. As $\left| {H:C} \right| < \infty $, we see that $akH$ is a finitely generated $kC$-module and, as ${I_1} = (kF{\bigcap I _1})k{K_1}$, we can conclude that $akH(kF{\bigcap I _1})kC = akH{I_1}$. Then the assertion follows from \cite[Lemma 2.4.1(ii)]{Tush2000}. 
\par  
	Step 2. Consider now the general case. The element $b$ may be written in the form $b = {a_1}{g_1} + ... + {a_n}{g_n}$ , where ${a_i} \in akL$ and ${g_i} \in N$. Then it is easy to note that $bkL \cong ( \oplus _{i = 1}^n{a_i}kL{g_i})/U$, where $U$ is a $kL$-submodule of $ \oplus _{i = 1}^n{a_i}{g_i}kL$, and hence $bkL/bkLI \cong ( \oplus _{i = 1}^n({a_i}kL/{a_i}kLI){g_i})/{U_1}$, where ${U_1}$ is a $kL$-submodule of $ \oplus _{i = 1}^n({a_i}kL/{a_i}kLI){g_i}$. Since $A$ is a central subgroup of $N/{I^\dag }$ , we see that $({a_i}kL/{a_i}kLI){g_i}) \cong {a_i}kL/{a_i}kLI$ considered as $kA$-modules and hence $bkL/bkLI \cong ( \oplus _{i = 1}^n({a_i}kL/{a_i}kLI))/{U_1}$, where ${U_1}$ is a $kL$-submodule of $ \oplus _{i = 1}^n({a_i}kL/{a_i}kLI)$. As it proved in step 1, ${a_i}kL/{a_i}kLI$ has a finite series each of whose quotient is isomorphic to some section of $akL/akLI$ and the assertion follows.
\par  
(iii). Step 1. At first, we show that $bkL$ contains a submodule $dkL$ isomorphic to $cgkL = ckLg$, where $0 \ne c \in akL$ and $g \in N$. The element $b$ may be written in the form $b = {a_1}{g_1} + ... + {a_n}{g_n}$ , where ${a_i} \in akL$ and ${g_i} \in N$. The proof is by induction on $n$. 
\par  
Evidently, $bkL + ({a_2}{g_2} + ...{a_n}{g_n})kL = {a_1}{g_1}kL + ({a_2}{g_2} + ...{a_n}{g_n})kL$ and hence $bkL + ({a_2}{g_2} + ...{a_n}{g_n})kL/({a_2}{g_2} + ...{a_n}{g_n})kL = {a_1}{g_1}kL + ({a_2}{g_2} + ...{a_n}{g_n})kL/({a_2}{g_2} + ...{a_n}{g_n})kL$. If $bkL \cap ({a_2}{g_2} + ...{a_n}{g_n})kL = 0$ then $bkL \cong {a_1}{g_1}kL/({a_1}{g_1}kL \cap ({a_2}{g_2} + ...{a_n}{g_n})kL)$. Then it follows from (i) that ${a_1}{g_1}kL \cap ({a_2}{g_2} + ...{a_n}{g_n})kL = 0$ and hence $bkL \cong {a_1}{g_1}kL$. So, we can put $c = {a_1}$ and $g = {g_1}$. 
\par  
If there is $0 \ne {b_1} \in bkL \cap ({a_2}{g_2} + ...{a_n}{g_n})kL$ then ${b_1}$ can be written in the form ${b_1} = {e_2}{g_2} + ...{e_n}{g_n}kL$ where ${e_i} \in akL$. As $0 \ne {b_1}kL \leqslant bkL$, changing  the element $b$ by ${b_1}$ we can apply the induction hypothesis. 
\par   
	 Step 2. By (ii), $dkL/dkLI$ has a finite series each of whose quotient is isomorphic to some section of $bkL/bkLI$ considered as $kA$-modules and it follows from step 1 that $ckLg/ckLgI$ has such a series. Since $A$ is a central subgroup of $N/{I^\dag }$, we can conclude that $ckLg/ckLgI$ and $ckL/ckLI$ are isomorphic as $kA$-modules. 
\end{proof}

Let $R$ be a ring, let $V$ be an $R$-module and $U$ be a submodule  of  $V$. An ideal $I$ of $R$ culls $U$ in $V$ if $VI < V$ and  $VI + U = V$ (see \cite{Broo88}). 


\begin{lemma}
Let $R$ be a ring, let $V$ be an $R$-module and $U$ be a submodule of $V$. An ideal $I$ of $R$ culls $U$ in $V$ if and only if ${I^n}$ culls $U$ in $V$ for any $n \in \mathbb{N}$. 
\end{lemma}

	\begin{proof} If $VI < V$ then $V{I^n} < V$ for any $n \in \mathbb{N}$. If $VI + U = V$ then $V{I^2} + UI = VI$  and hence $V{I^2} + U = V{I^2} + UI + U = VI + U = V$. Thus, $V{I^2} + U = V$ and by induction we can conclude that $V{I^n} + U = V$ for any $n \in \mathbb{N}$. 
\end{proof}

	Let $R$ be a domain, $K$ be a normal subgroup of a group $H$. According to \cite[Introduction, p. 89]{Broo88}, we say that that an $RK$-module $0 \ne V$ is $H$-ideal critical if for any submodule $0 \ne {V_1} \leqslant V$ there is an $H$-invariant ideal $I$ of $RK$ such that: 

\begin{romanlist}
\item $I$ culls ${V_1}$ in $V$;
\item $I$ has the weak Artin-Rees property that is, if $U$ is any finitely generated $RK$-module with submodule ${U_1}$ then $U{I^n} \cap {U_1} \leqslant {U_1}I$ for some $n \in \mathbb{N}$ dependent on $U$ and ${U_1}$. 
\end{romanlist}


\begin{proposition}
Let $N$ be a minimax nilpotent torsion-free normal subgroup of a soluble group $G$ such that ${r_0}(G) < \infty $ and let $K$ be a $G$-invariant subgroup of $N$.   Let $H$ be a finitely generated subgroup of $K$, let $k$ be a field of characteristic zero and let  $0 \ne J$ be a right ideal of $kH$. Then:     
\begin{romanlist}

\item there is a $G$-large $N$-central  ideal $I$ of $kK$ such that $kH \cap I$ culls $J$ in  $kH$; 
\item if the group  $G = N$ then the ideal $I$ may be chosen such that  $I = (I \cap kF)kK$, where $K \geqslant F \geqslant {I^\dag }$,  $F/{I^\dag }$ is a central section of  $G$ and  $I \cap kF$ is a maximal ideal of $kF$
\item if  $H = N$ then $kK$ is $H$-ideal critical uniform $kK$-module. 
\end{romanlist}
\end{proposition}
\begin{proof} (i), (ii). By \cite[Theorem 3.7]{Wehr09}, the ring $kH$ is Noetherian and hence the ideal $J$ is finitely generated. Therefore, there is a finitely generated subfield ${k_1} \leqslant k$ such that $J = (J \cap {k_1}H)kH$. By \cite[Proposition 2.4.3]{Tush2000}, there is a $G$-large $N$-central ideal ${I_1}$ of ${k_1}K$ such that ${k_1}H \cap {I_1}$ culls $J \cap {k_1}H$ in ${k_1}H$ and if $G = N$ then the ideal $I_1$ may be chosen such that ${I_1} = ({I_1} \cap {k_1}F){k_1}K$, where $K \geqslant F \geqslant I_1^\dag $, $F/I_1^\dag $ is a central section of $G$ and ${I_1} \cap {k_1}F$ is a maximal ideal of ${k_1}F$. Put $I = ({I_1} \cap {k_1}F)kK$ then $kH \cap I$ culls $J$ in  $kH$. If $G = N$ and $I \cap kF$ is not  a maximal ideal of $kF$ then we can change $I \cap kF$ by any maximal ideal $P$ of  $kF$ containing $I \cap kF$ and put $I = PkK$.  
\par  
(iii). By \cite[Lemma 2.4.1(i)]{Tush2000}, any $H$-large ideal $I$ of $kK$ is polycentral and it follows from \cite[Theorem 6.12]{Wehr09} that $I$ has the weak Artin-Rees property. Then (i) implies that $kK$ is $H$-ideal critical. It follows from \cite[Corollary 37.11]{Pass89} that $kK$ is an Ore domain and hence $kK$ is uniform. 
\end{proof}

Let $K$ be a normal subgroup of a  group $N$  such that the quotient group $N/K$ is torsion-free abelian of finite rank. Let ${I^\dag }$ be an $N$-invariant subgroup of $K$ such that the quotient group $K/{I^\dag }$ is finite. As the quotient group $K/{I^\dag }$ is finite and the quotient group $N/K$ is torsion-free abelian of finite rank, it follows from Lemma 2.2(i) that $N/{I^\dag }$ has a torsion-free abelian characteristic central subgroup $A$ of finite index. For any subgroup $X$ of $N$ such that $K \leqslant X$ we denote  ${A_X} = A \cap X/{I^\dag }$.  


\begin{lemma}
Let $K$ be a normal subgroup of a torsion-free finitely generated nilpotent  group $H$ such that the quotient group $H/K$ is free abelian and let $\chi $ a full system of subgroups of $H$ over $K$. Let  $k$ be a field and let $0 \ne W$ be a $kH$-module which is $kK$-torsion-free. Then :
\begin{romanlist}
\item for any element $0 \ne a \in W$ there is a  right ideal ${J_0}$ of $kK$  such that if $I$ is an $H$-large ideal of $kK$ which culls ${J_0}$ in $kK$ then $X \in \chi (akH)$ if and only if  the quotient module $akH/akHI$ is not $k{A_X}$-torsion, where $A$ is a characteristic central torsion-free subgroup of finite index in $H/{I^\dag }$ and ${A_X} = A \cap X/{I^\dag }$; 
\item  for any element $0 \ne a \in W$ any subgroup $X \in \chi (akH)$ there is a  right ideal ${J_X}$ of $kK$  such that if $I$ is an $H$-large ideal of $kK$ which culls ${J_X}$  in $kK$ then for any element $0 \ne c \in akX$ the quotient module $ckH/ckHI$ is not $k{A_X}$-torsion, where $A$ is a characteristic central torsion-free subgroup of finite index in $H/{I^\dag }$ and ${A_X} = A \cap X/{I^\dag }$; 
\item there are a  cyclic $kH$-submodule $0 \ne akH \leqslant W$ and a right ideal $J$ of $kK$  such that if $I$ is an $H$-large ideal of $kK$ which culls $J$ in $kH$ then  for any cyclic $kH$-submodule  $0 \ne bkH \leqslant akH$ we have $X \in \chi (bkH)$ if and only if  the quotient module $bkH/bkHI$ is not $k{A_X}$-torsion, where $A$ is a characteristic central torsion-free subgroup of finite index in $H/{I^\dag }$ and ${A_X} = A \cap X/{I^\dag }$. 
\end{romanlist}
\end{lemma}
\begin{proof} (i) Since the module $0 \ne W$ is $kK$-torsion-free, $V = akK \cong k{K_{kK}}$ is an $H$-ideal critical $kK$-module by Proposition 3.1(iii). 
\par  
It follows from \cite[Lemma 8]{Broo88} that there is a right ideal ${J_0}$ of $kK$ such that if an $H$-large ideal $I$ of $kK$ culls ${J_0}$ in $kK$ (i.e. $I$ culls $a{J_0}$ in $V = akK \cong k{K_{kK}}$)  then $X \in \chi (bkH)$ if and only if  $bkH/bkHI{ \otimes _{k{A_X}}}{Q_X} \ne 0$, where ${Q_X}$ denotes the field of fractions of the domain $k{A_X}$. It is well known that the relation $(bkH/bkHI){ \otimes _{k{A_X}}}{Q_X} \ne 0$ means that the quotient module $bkH/bkHI$ is not $k{A_X}$-torsion. 
\par  
(ii) Let $X \in \chi (akH)$, arguments of the proof of \cite[Lemma 14]{Broo88} show that there is a right ideal ${J_X} \leqslant kK$ such that if an $H$-large ideal $I$ of $kK$ culls ${J_X}$ in $kK$ then $akXI = akHI \cap akX$. By Lemma 3.4, ${I^n}$ also culls ${J_X}$ and we can conclude that 
\begin{equation}
                            akX{I^n} = akH{I^n} \cap akX.                      \label{3.2}
\end{equation}
    It follows from the definition of $M\chi (akH)$ that $akX \cong kX$ and it easily follows from Proposition 3.1(iii) that ${ \cap _{n \in N}}akX{I^n} = 0$. Then there is $m \in \mathbb{N}$ such that $c \in akX{I^{m - 1}}\backslash akX{I^m}$. Therefore, $ckX/(ckX \cap akX{I^m}) \cong (ckX + akX{I^m})/akX{I^m}$ is a non-zero submodule of $akX{I^{m - 1}}/akX{I^m}$. 
\par  
     Show now that $S  =  kX{I^{m - 1}}/kX{I^m}$ is a torsion-free $k{A_X}$-module and hence so is $akX{I^{m - 1}}/akX{I^m}$. Evidently, $S = ({I^{m - 1}}/{I^m}){ \otimes _{kC}}kB$, where $C = K/{I^\dag }$ and $B = X/{I^\dag }$. Then, as $A \cap C  =  1$,  $S$ is $k{A_X}$-torsion-free. 
\par  
     Thus, $akX{I^{m - 1}}/akX{I^m}$ is a $k{A_X}$-torsion-free module and hence so is $(ckX + akX{I^m})/akX{I^m} \cong ckX/(ckX \cap akX{I^m})$. Due to (\ref{3.2}), as $c \in akX$, we see that $akX{I^n} \cap ckX = akH{I^n} \cap ckX$ and hence $ckX/(ckX \cap akH{I^m})$ is a $k{A_X}$-torsion-free module. As $ckX \leqslant akX{I^{m - 1}}$ , we see that $ckH \leqslant akH{I^{m - 1}}$ and hence  $ckHI \leqslant akH{I^m}$. It implies that $ckX \cap ckHI \leqslant ckX \cap akH{I^m}$ and , as $ckX/(ckX \cap akH{I^m})$ is a $k{A_X}$-torsion-free module, the  quotient module $ckX/(ckX \cap ckHI) \simeq (ckX + ckH)/ckHI \leqslant ckH/ckHI$ is not $k{A_X}$-torsion. 
\par  
(iii). We apply Lemma 3.2 to the module $W$ in the case, where $L = H$. Then we can choose an element $0 \ne a \in W$ such that for any $0 \ne b \in akH$ the following relations hold:
\begin{equation}
                                   \chi (akH) = \chi (bkH)                                 \label{3.3}
\end{equation}  
and
\begin{equation}
                                        akX \cap bkH \ne 0                                  \label{3.4}
\end{equation}
for any $X \in M\chi (akH)$.     
\par  
	 Let ${J_0}$ be a right ideal of $kK$ from (i) defined for the element $a$ and let ${J_X}$ be a right ideal  of $kK$ from (ii) defined 
for the element $a$ and for a subgroup $X \in \chi (akH)$. Put $J = {J_0} \cap ({ \cap _{X \in \chi (akH)}}{J_X})$ and let $I$ be an 
$H$-large ideal of $kK$ which culls $J$ in $kK$ then, by \cite[Lemma 6]{Broo88}, $I$ culls ${J_0}$ and ${J_X}$ in $kK$ for each 
$X \in \chi (akH)$. Therefore, (i) and (ii) hold for the module $akH$ and the ideal $I$. 
\par  
	Let $0 \ne bkH \leqslant akH$ and suppose that for some $X \in \chi $ the quotient module $bkH/bkHI$ is not $k{A_X}$-torsion. Then it follows from Lemma 3.3(ii) that the quotient module $akH/akHI$ is not $k{A_X}$-torsion and by (i) $X \in \chi (akH)$. Therefore, by (\ref{3.3}), $X \in \chi (bkH)$. 
\par  
	Let now $X \in \chi (bkH)$ and show that the quotient module $bkH/bkHI$ is not $k{A_X}$-torsion. Evidently, $X \leqslant Y$ for some $Y \in M\chi (bkH)$ and if we show that $bkH/bkHI$ is not $k{A_Y}$-torsion then $bkH/bkHI$ is not $k{A_X}$-torsion because ${A_X} \leqslant {A_Y}$. Thus, we can assume that $X \in M\chi (bkH)$ and hence, by (\ref{3.4}), there is $0 \ne c \in akX \cap bkH$. As $0 \ne c \in akX$, it follows from (ii) that the quotient module $ckH/ckHI$ is not $k{A_X}$-torsion. Then, as $0 \ne ckH \leqslant bkH$, it follows from Lemma 3.3(ii) that the quotient module $bkH/bkHI$ is not $k{A_X}$-torsion.  
\end{proof}


\begin{proposition}
Let $K$ be a normal subgroup of a torsion-free minimax nilpotent  group $N$ such that the quotient group $N/K$ is torsion-free abelian. Let $L$ be a dense subgroup of $N$ such that $K \leqslant L$ and  the quotient group $L/K$ is free abelian.  Let  $k$ be a field and let $0 \ne W$ be a  $kN$-module which is $kK$-torsion-free. Then there are an element $0 \ne a \in W$, a finitely generated dense subgroup $H$ of $L$ and a right ideal $J$ of $k(H \cap K)$  such that if $I$ is an $N$-large ideal of $kK$ such that $I \cap k(K \cap H)$ culls $J$ in $k(H \cap K)$ then  for any $0 \ne b \in akN$ we have $X \in \chi (bkL)$ if and only if  the quotient module $bkL/bkLI$ is not $k{A_X}$-torsion, where ${A_X} = A \cap X/{I^\dag }$ and $A$  is a characteristic central torsion-free subgroup of finite index in $N/{I^\dag }$. 
\end{proposition}

\begin{proof} By Lemmas 3.2 and 3.3, we can choose a cyclic submodule $0 \ne akN \leqslant W$ which satisfies the conditions (i)-(iv) of Lemma 3.2 and the conditions (ii),(iii) of Lemma 3.3  and these conditions hold for any cyclic submodule $0 \ne bkN \leqslant akN$. 
\par  
Step 1. Consider at first the case where $N = L$. By Lemma 3.2(i), there is a finitely generated dense subgroup $H$ of $L$ such that
\begin{equation}
 akL = akH{ \otimes _{kH}}kL.                                                 \label{3.5}  
\end{equation} 
It follows from Lemma 3.1 that the subgroup $H$ may be chosen such that for any $L$-invariant subgroup ${I^\dag } \leqslant K$ of finite index in $K$
\begin{equation}
                                       {I^\dag }H = L                           \label{3.6}  
\end{equation}
and
\begin{equation}
   \chi = \left\{ {{I^\dag }X|X \in {\chi _H}} \right\},                           \label{3.7} 
  \end{equation}
where ${\chi _H} = \left\{ {\left\langle {{K_1},\{ {x_j}|j \in J\} } \right\rangle | J \subseteq \left\{ {1,..., n} \right\}} \right\}$ is a full system of subgroups of $H$ over ${K_1} = K \cap H$.    
\par  
Show that 
\begin{equation}
   \chi (akL) = \left\{ {{I^\dag }X|X \in {\chi _H}(akH)} \right\}.     \label{3.8}
\end{equation}

We use Eq. (\ref{3.7}). Let $X \in {\chi _H}(akH)$ if ${I^\dag }X \in \chi \backslash \chi (akL)$ then, by Lemma 3.2(ii), the module $akL$ is $k({I^\dag }X)$-torsion and, as $X$ is a dense subgroup of ${I^\dag }X$, it follows from \cite[Lemma 2.2.3(ii)]{Tush2000} that the module $akL$ is $kX$-torsion. But then $X \notin {\chi _H}(akH)$ and a contradiction is obtained. Suppose now that ${I^\dag }X \in \chi (akL)$, where $X \in {\chi _H}$ then, by Lemma 2.3(ii),  the module $akL$ is $k({I^\dag }X)$-torsion-free. Therefore, $akH$ is $kX$-torsion-free and hence $X \in {\chi _H}(akH)$. So, the relation (\ref{3.8}) holds. 
\par  
By Lemma 3.5(iii), we can choose the element $a$ for which there is a right ideal $J$ of $k{K_1}$  such that if ${I_1}$ is an $H$-large ideal of $k{K_{}}$ which culls $J$ in $k{K_1}$ then for any cyclic $kH$-submodule $0 \ne bkH \leqslant akH$ we have $X \in {\chi _H}(bkH)$ if and only if the quotient module $bkH/bkH{I_1}$ is not $k{A_X}$-torsion, where ${A_X} = A \cap X/{I_1}^\dag $ and $A$ is a central subgroup of finite index in $H/{I^\dag }$. 
\par  
	Let $I$ be an $L$-large ideal of $kK$ such that ${I_1} = I \cap k{K_1}$ culls $J$ in $k{K_1}$. It follows from (\ref{3.6}) and \cite[Lemma 2.2.6]{Tush2000} that $akL/akLI$ and $akH/akH{I_1}$ are isomorphic as $kA$-modules, where $A$ is a central subgroup of finite index in $L/{I^\dag }$. Then it follows from Lemma 3.5(iii) and (\ref{3.8}) that $X \in \chi (akL)$ if and only if the quotient module $akL/akLI$ is not $k{A_X}$-torsion, where ${A_X} = A \cap X/{I^\dag }$ and $A$  is a central subgroup of finite index in $L/{I^\dag }$. 
\par  
 	Let $0 \ne b \in akL$, by Lemma 3.2(iv), there is an element $0 \ne c \in bkL \cap akH$ and it follows from (\ref{3.5}) that we have a chain of submodules $ckH{ \otimes _{kH}}kL = ckL \leqslant bkL \leqslant akL = akH{ \otimes _{kH}}kL$. By Lemma 3.2(iii), $\chi (ckL) = \chi (bkL) = \chi (akL)$. Since $ckL = ckH{ \otimes _{kH}}kL$, the same arguments as in the case $akL$ show that $X \in \chi (ckL) = \chi (bkL) = \chi (akL)$  if and only if $ckL/ckLI$ is not $k{A_X}$-torsion. Let $X \in \chi (ckL) = \chi (bkL) = \chi (akL)$ then $ckL/ckLI$ is not $k{A_X}$-torsion, and it follows from Lemma 3.3 that $bkL/bkLI$ is not $k{A_X}$-torsion because $ckL \leqslant bklL$. Suppose now that $bkL/bkLI$ is not $k{A_X}$-torsion then it follows from Lemma 3.3 that $akL/akLI$ is not $k{A_X}$-torsion because $bkL \leqslant akl$ and hence $X \in \chi (ckL) = \chi (bkL) = \chi (akL)$.
\par  
	Step 2. Consider now the general case. By step 1, we can chose the element $0 \ne a \in W$ such that the assertion holds for any element $0 \ne b \in akL$. Besides, as it was mention in the beginning of the proof, $akN$ satisfies the conditions (i)-(iv) of Lemma 3.2 and the conditions (ii), (iii) of Lemma 3.3 and these conditions hold for any cyclic submodule $0 \ne bkN \leqslant akN$. 
\par  
	Let $X \in \chi (bkL)$ and suppose that the quotient module $bkL/bkLI$ is $k{A_X}$-torsion. By Lemma 3.3(iii), there is a $kL$-submodule $0 \ne ckL \leqslant akL$ such that $ckL/ckLI$ has a finite series each of whose quotient is isomorphic to some section of $bkL/bkLI$ considered as $kA$-modules. Therefore, the quotient module $ckL/ckLI$ is $k{A_X}$-torsion. On the other hand, by Lemma 3.2(iii), $X \in \chi (ckL)$ and by, step 1, the quotient module $ckL/ckLI$ is not $k{A_X}$-torsion. Thus, a contradiction is obtained and hence the quotient module $bkL/bkLI$ is not $k{A_X}$-torsion. 
\par  
	Suppose now that the quotient module $bkL/bkLI$ is not $k{A_X}$-torsion. By Lemma 3.3(ii), $bkL/bkLI$ has a finite series each of whose quotient is isomorphic to some section of $akL/akLI$ considered as $kA$-modules and hence the quotient module $akL/akLI$ is not $k{A_X}$-torsion. Then, by step 1, $X \in \chi (akL)$ and, by Lemma 3.2(iii),  $X \in \chi (bkL)$. 
\end{proof}


\section{A Finite Set of Commutative Invariants for Modules over Group Rings of Nilpotent Minimax Groups}

Let $R$ be a commutative Noetherian ring and let $I$ be an ideal of $R$. Let ${\mu _R}(I)$ be the set of prime ideals of $R$ minimal over $I$, by \cite[Chap. II, \S 4, Corollary 3]{Bour}, the set ${\mu_R}(I)$ is finite. 
\par  
Let $P$ be a prime ideal of a commutative ring $R$ and let ${R_P}$ be the localization  of $R$ at the ideal $P$. Let $M$ be an $R$-module, the support $Sup{p_R}M$ of the module $M$ consists of all prime ideals $P$ of $R$ such that ${M_P} = M{ \otimes _R}{R_P} \ne 0$ (see \cite[Chap. II, \S 4]{Bour}). By \cite[Chap. IV, \S 1, Theorem 2]{Bour}, if $R$ and $M$ are Noetherian then the set ${\mu_R}(M)$  of minimal elements of $Sup{p_R}M$ coincides with ${\mu_R}(Ann_R(M))$, where ${\mu_R}(Ann_R (M))$ is the set of prime ideals of $R$ which are minimal over $An{n_R}M$. Thus, we have
\begin{equation}
{\mu_R}(M) = {\mu_R}(Ann_{R}(M)).      \label{4.1} 
\end{equation}


\begin{lemma}
Let $A$ be a finitely generated abelian group, let $k$ be a field and let  $I$ be an ideal of $kA$. Then for any  subgroup $B\leq A$ of finite index in $A$ 
the following assertions hold: 
\begin{romanlist} 
\item for any ${Q \in \mu _{kB}}(I \cap kB)$ there is a prime ideal $P$ of $kA$  such  that  $I\leq P$, $Q = P \cap kB$ and any such 
prime  $P$ is necessarily in  ${\mu _{kA}}(I)$;
\item if $P \in  {\mu _{kA}}(I )$ then $P \cap kB \notin  {\mu _{kB}}(I \cap kB )$ if and only if there is $S \in  {\mu _{kA}}(I)$ such that 
$S \cap kB < P \cap kB$;
\item if $J$ is an ideal of $kA$ such that $\mu _{kA}(I ) = \mu _{kA}(J )$ then $\mu _{kB}(I \cap kB) = \mu _{kB}(J \cap kB)$;
\item if $M$ and $M_1$ are finitely generated $kA$-modules such that $\mu _{kA}(M) = \mu _{kA}(M_1)$ then $\mu _{kB}(M) = \mu _{kB}(M_1)$.
\end{romanlist}
\end{lemma}

\begin{proof} (i) Evidently, the quotient ring $kA/I$ is integer over $(kB+I)/I \cong kB/(kB \cap I)$. Then it follows from \cite[Ch.V, \S 2, Theorem 3]{ZaSa58} that there is a prime ideal $P \leq kA$ such that  $I \leq P$  and $P \cap kB = Q$. If $P \notin  {\mu _{kA}}(I)$ then there is an ideal $ P' \in  {\mu _{kA}}(I )$ such that $I \leq P' < P$ and hence   $I \cap kB \leq P'  \cap kB  \leq P  \cap kB =Q$ and, as   ${Q \in \mu _{kB}}(I \cap kB) $, we can conclude that $I \cap kB \leq P'  \cap kB  = P  \cap kB =Q$. Since $I \leq P' < P$, it contradicts Complement 1 to \cite[Ch.V, \S 2, Theorem 3]{ZaSa58}. Thus,  $P \in  {\mu _{kA}}(I)$. 
\par 
 (ii) Let  $P \in  {\mu _{kA}}(I )$. Suppose that  $P \cap kB \notin  {\mu _{kB}}(I \cap kB )$ then there is a prime ideal $Q \in  {\mu _{kB}}(I \cap kB )$ such that $Q < P \cap kB$. By (i), there is $S \in  {\mu _{kA}}(I )$ such that  $Q = S \cap kB$ and hence $S \cap kB < P \cap kB$. Suppose now that there 
 is $S \in  {\mu _{kA}}(I)$ such that $S \cap kB < P \cap kB$. Then $I \cap kB \leq S \cap kB < P \cap kB$ and hence $P \cap kB \notin  {\mu _{kB}}(I \cap kB )$. 
 \par  
 (iii) Let $Q \in \mu _{kB}(I \cap kB)$. By (i), there is $P \in  {\mu _{kA}}(I )$ such that  $Q = P \cap kB$. Then, by (ii), there are no 
 $S \in  {\mu _{kA}}(I)$ such that $S \cap kB < P \cap kB$. Since $\mu _{kA}(I ) = \mu _{kA}(J )$, we can conclude that there is $P \in  {\mu _{kA}}(J)$ 
 such that $Q = P \cap kB$ and there are no $S \in  {\mu _{kA}}(J)$ such that $S \cap kB < P \cap kB$. Then it follows from (ii) that 
 $Q \in \mu _{kB}(J \cap kB)$. The same arguments show that if $Q \in \mu _{kB}(J \cap kB)$ then $Q \in \mu _{kB}(I \cap kB)$ and the assertion follows. 
 \par
 (iv) It follows from (\ref{4.1}) that  $\mu_{kA}(Ann_{kA}(M))=\mu_{kA}(M) = \mu_{kA}(M_1) = \mu_{kA}(Ann_{kA}(M_1))$. Then, as 
 $\mu_{kA}(Ann_{kA}(M))= \mu_{kA}(Ann_{kA}(M_1))$, it follows from (iii) that $\mu_{kB}(Ann_{kA}(M)\cap kB) = {\mu_{kB}}(Ann_{kA}(M_1)\cap kB)$. 
 By (\ref{4.1}), we have $\mu_{kB}(M) = \mu_{kB}(Ann_{kB}(M))= \mu_{kB}(Ann_{kA}(M)\cap kB)$ 
 and $\mu_{kB}(M_1) = \mu_{kB}(Ann_{kB}(M_1)) = {\mu_{kB}}(Ann_{kA}(M_1)\cap kB)$. Then, as  $\mu_{kB}(Ann_{kA}(M)\cap kB) = {\mu_{kB}}(Ann_{kA}(M_1)\cap kB)$, 
 we can conclude that  $\mu_{kB}(M) = \mu_{kB}(M_1)$. 
  \end{proof}

	Let $G$ be a group and let $K$ be a normal subgroup of $G$. Let $L$ be a subgroup of $G$ such that $K \leqslant L \leqslant G$ and the quotient group $L/K$ is finitely generated free abelian. Let $k$ be a field and let $I$ be a $G$-large ideal of the group ring $kK$. Then the quotient group $K/I^\dag$ is finite and, as the quotient group $L/K$ is finitely generated free abelian, it follows from Lemma 2.2(i) that
	 $L/I^\dag$  has a characteristic central torsion-free subgroup $A$ of finite index. Since the quotient group $L/K$ is finitely generated free abelian, the subgroup $A$ is finitely generated free abelian. By \cite[Theorem 3.7]{Wehr09}, the group ring $kA$ is a Noetherian. Let $W$ be a finitely generated $RL$-module then $W/WI$ is a finitely generated $kA$ -module and hence $W/WI$ is a Noetherian $kA$-module. So, the finite set ${\mu _{kA}}(W/WI) = {\mu _{kA}}(An{n_{kA}}(W/WI))$ is defined. 
\par  
If $P$ is a prime ideal of the group ring $kA$ then the dimension $d(P)$ of $P$ is the Krull dimension of the quotient ring $kA/P$. It is well known that the dimension $d(P)$ coincides with $r(X)$, where $X$ is a maximal subgroup of $A$ such that $kX \cap P = 0$. Thus, the dimension of any prime ideal of $kA$ is bounded with $r(A)$ and $d(\mu)$ will denote the maximal dimension of ideals from $\mu$, where $\mu$ is a set  of prime ideals of $kA$ (if $\mu = \emptyset $ we put $d(\mu) = -1$). 


\begin{lemma}
Let $L$ be a minimax torsion-free nilpotent group and let $K$ be a normal subgroup of $L$ such that the quotient group $L/K$ is free abelian.  
Let $T$ be a subgroup of finite index in  $L$ such that $K \leqslant T$. Let $k$ be a field and let $W$be a cyclic $kL$-module which is $kK$-torsion-free. Let $I$ be an $L$-large ideal of $kK$ such that $W \ne WI$. Let $A$ be a central torsion-free subgroup of finite index in $L/{I^\dag }$ and let $B$ be a subgroup of finite index in $A \cap (T/{I^\dag })$ .  Then:
\begin{romanlist}

\item $Sup{p_{kB}}  {W_1}/{W_1}I \subseteq Sup{p_{kB}}  W/WI$  for any cyclic $kT$-submodule  ${W_1} \leqslant W$; 
\item  for any  cyclic $kL$-submodule ${W_1}$ of  $W$ and any $P \in {\mu _{kA}}({W_1}/{W_1}I)$ there is  $Q \in {\mu _{kA}}(W/WI)$ such that $Q \leqslant P$; 
\item ${\mu _{kA}}({W_2}/{W_2}I) \cap {\mu _{kA}}(W/WI) \subseteq {\mu _{kA}}({W_1}/{W_1}I) \cap {\mu _{kA}}(W/WI)$ for any cyclic $kL$-submodules  ${W_1}$ and ${W_2}$ of $W$ such that ${W_2} \leqslant {W_1}$; 
\item there exists a cyclic $kL$-submodule $0 \ne V \leqslant W$ such that  ${\mu _{kB}}({V_1}/{V_1}I) = {\mu _{kB}}(V/VI)$  for any cyclic $kL$-submodule  $0 \ne {V_1} \leqslant V$ and any  subgroup $B\leq A$ of finite index in $A$. 
\end{romanlist}
\end{lemma}

     \begin{proof} (i). Since $\left| {L:T} \right| < \infty $, we see that $akL$ is a finitely generated $kT$-module. If the group $L$ is finitely generated, the assertion is proved in \cite[Lemma 5(i)]{Tush02}. 
\par  
     Consider now the general case. Suppose that $W = akL$ and ${W_1} = bkT$, where $b \in akL$. By \cite[Lemma 2.2.5]{Tush2000}, there is a finitely generated subgroup $H \leqslant L$ such that $W = akH{ \otimes _{kH}}kL$. Evidently, taking the subgroup $H$ bigger if it is necessary, we can assume that $b \in akH$ and $L = {I^\dag }(K \cap H)$. Therefore, ${W_1} = bk(H \cap T){ \otimes _{k(H \cap T)}}kT$ and, as $L = {I^\dag }(K \cap H)$, we have $L = {I^\dag }H$ and $T = {I^\dag }(H \cap T)$ . Then it follows from   \cite[Lemma 2.2.6]{Tush2000}   that there are $kH$-module isomorphism $akL/akLI \simeq akH/akHI'$ and $k(H \cap T)$-module isomorphism  $bkT/bkTI \simeq bk(H \cap T)/bk(H \cap T)I'$, where $I' = kH \cap I$ and these isomorphisms induce $kB$-module isomorphisms. So, the assertion follows from the considered above case, where the group $L$ is finitely generated. 
\par  
    (ii). We apply (i) in the case where $T = L$ and $A = B$  then $Sup{p_{kA}}{W_1}/{W_1}I \subseteq Sup{p_{kA}}W/WI$ and hence ${\mu _{kA}}({W_1}/{W_1}I) \subseteq Sup{p_{kA}}{W_1}/{W_1}I \subseteq Sup{p_{kA}}W/WI$. Then the assertion follows from the definition of ${\mu _{kA}}(W/WI)$.      
\par  
(iii). By (ii), for any $P \in {\mu _{kA}}({W_2}/{W_2}I) \cap {\mu _{kA}}(W/WI)$ there is  $Q \in {\mu _{kA}}({W_1}/{W_1}I)$ such that $Q \leqslant P$ and, also by (ii), there is ${P_1} \in {\mu _{kA}}(W/WI)$ such that ${P_1} \leqslant Q \leqslant P$. However, $P \in {\mu _{kA}}(W/WI)$ and it follows from the minimality of $P$ that ${P_1} = Q = P$. Thus, $P =   Q \in {\mu _{kA}}({W_1}/{W_1}I)$ and hence $P \in {\mu _{kA}}({W_1}/{W_1}I) \cap {\mu _{kA}}(W/WI)$. 
\par  
(iv)To simplify denotations, we denote $\mu_{kA}$ my $\mu$.  By  Lemma 4.1(iv), it is sufficient to show that there exists a cyclic $kL$-submodule $0 \ne V \leqslant W$ such that  
${\mu }({V_1}/{V_1}I) = {\mu }(V/VI)$  for any cyclic $kL$-submodule  $0 \ne {V_1} \leqslant V$.
\par
Suppose that the assertion is not true then there is a descending chain $\{ W_n|n\in\mathbb{N}\}$ of cyclic submodules of $W$ such that 
$${\mu }({W_n}/{W_n}I) \ne {\mu }({W_{n + 1}}/{W_{n + 1}}I)$$ 
for each $n \in \mathbb{N}$. 
\par  
Step 1. Note at first that ${\mu }({W_t}/{W_t}I) \ne {\mu }({W_n}/{W_n}I)$ if $t \ne n$. Suppose that there are $t>n$ such that $\mu
(W_{t}/W_{t})=\mu (W_{n}/W_{n}I)$, we can assume that $t$ and $n$ are
nearest numbers with this property.
Then it follows from the definition of the chain $\{ {W_n}|n \in \mathbb{N}\} $ that $t > n + 1 > n$ and ${\mu }({W_t}/{W_t}I) \ne {\mu }({W_{n + 1}}/{W_{n + 1}}I) \ne {\mu }({W_n}/{W_n}I)$. Since ${W_t} \leqslant {W_{n + 1}} \leqslant {W_n}$, it follows from (iii) that ${\mu }({W_t}/{W_t}I) \cap {\mu }({W_n}/{W_n}I) \subseteq {\mu }({W_{n + 1}}/{W_{n + 1}}I) \cap {\mu }({W_n}/{W_n}I)$ and, as ${\mu }({W_t}/{W_t}I) = {\mu }({W_n}/{W_n}I)$, we have ${\mu }({W_n}/{W_n}I) \subseteq {\mu }({W_{n + 1}}/{W_{n + 1}}I) \cap {\mu }({W_n}/{W_n}I)$. Therefore, ${\mu }({W_n}/{W_n}I) \subseteq {\mu }({W_{n + 1}}/{W_{n + 1}}I)$. By (ii), for any $P \in {\mu }({W_{n + 1}}/{W_{n + 1}}I)$ there is  $Q \in {\mu }({W_n}/{W_n}I)$ such that $Q \leqslant P$ and, as ${\mu }({W_n}/{W_n}I) \subseteq {\mu }({W_{n + 1}}/{W_{n + 1}}I)$, it follows from the minimality of $P$ that $Q = P$. It implies that ${\mu }({W_n}/{W_n}I) \supseteq {\mu }({W_{n + 1}}/{W_{n + 1}}I)$ and hence ${\mu }({W_n}/{W_n}I) = {\mu }({W_{n + 1}}/{W_{n + 1}}I)$ but this contradicts the definition of the chain $\{ {W_n}|n \in N\} $. Thus, ${\mu }({W_t}/{W_t}I) \ne {\mu }({W_n}/{W_n}I)$ if $t \ne n$.
\par  
Step 2. It follows from (iii) that ${\mu }({W_t}/{W_t}I) \cap {\mu }({W_1}/{W_1}I) \subseteq {\mu }({W_m}/{W_m}I) \cap {\mu }({W_1}/{W_1}I)$ for any $t > m > 1$. Then, as the set ${\mu }({W_1}/{W_1}I)$ is finite, there is  $m \in \mathbb{N}$ such that 

$${\mu }({W_t}/{W_t}I) \cap {\mu }({W_1}/{W_1}I) = {\mu }({W_m}/{W_m}I) \cap {\mu }({W_1}/{W_1}I)$$               

for all $t > m$. We put ${W_m} = {U_2}$ and ${W_1} = {U_1}$. The same arguments show that there is  $m \in \mathbb{N}$ such that 

$${\mu }({W_t}/{W_t}I) \cap {\mu }({U_2}/{U_2}I) = {\mu }({W_m}/{W_m}I) \cap {\mu }({U_2}/{U_2}I)$$               

for all $t > m$. We put ${W_m} = {U_3}$. Continuing this process we obtain a descending chain  $\{ {U_n}|n \in \mathbb{N}\}  \subseteq \{ {W_n}|n \in \mathbb{N}\} $ such that 
\begin{equation}
      {\mu }({U_t}/{U_t}I) \cap {\mu }({U_n}/{U_n}I) = {\mu }({U_m}/{U_m}I) \cap {\mu }({U_n}/{U_n}I)     \label{4.2} 
\end{equation}
for all $t,m,n \in \mathbb{N}$ such that $t \geqslant m \geqslant n$. Since  ${\mu }({W_t}/{W_t}I) \ne {\mu }({W_n}/{W_n}I)$ if $t \ne n$, we see that  
\begin{equation}
  {\mu }({U_t}/{U_t}I) \ne {\mu }({U_n}/{U_n}I)                   \label{4.3}
\end{equation}
if $t \ne n$. 
\par
Show now that if 
$$d({\mu }({U_{n + 1}}/{U_{n + 1}}I)\backslash {\mu }({U_n}/{U_n}I)) \geqslant 0$$ 
(i.e. ${\mu }({U_{n + 1}}/{U_{n + 1}}I)\backslash {\mu }({U_n}/{U_n}I) \ne \emptyset $) then 
\begin{equation}
          d({\mu }({U_{n + 1}}/{U_{n + 1}}I)\backslash {\mu }({U_n}/{U_n}I)) < d({\mu }({U_n}/{U_n}I)\backslash {\mu }({U_{n - 1}}/{U_{n - 1}}I))      \label{4.4}
\end{equation}
for any $n \geqslant 2$. By (ii), for any $P \in {\mu }({U_{n + 1}}/{U_{n + 1}}I)\backslash {\mu }({U_n}/{U_n}I)$ 
there is  $Q \in {\mu }({U_n}/{U_n}I)$ such that $Q \leqslant P$ and, as $P \notin {\mu }({U_n}/{U_n}I)$, we have $Q < P$.  
By (\ref{4.2}), ${\mu }({U_{n + 1}}/{U_{n + 1}}I) \cap {\mu }({U_{n - 1}}/{U_{n - 1}}I) = {\mu }({U_n}/{U_n}I) \cap {\mu }({U_{n - 1}}/{U_{n - 1}}I)$ and hence if $Q \in {\mu }({U_{n - 1}}/{U_{n - 1}}I)$ then $Q \in {\mu }({U_{n + 1}}/{U_{n + 1}}I)$ 
but $P \in {\mu }({U_{n + 1}}/{U_{n + 1}}I)$ and it contradicts $Q < P$. Thus, $Q \notin {\mu }({U_{n - 1}}/{U_{n - 1}}I)$ 
and hence $Q \in \mu({U_n}/{U_n}I)\backslash \mu({U_{n - 1}}/{U_{n - 1}}I)$. 
So, for any $P \in {\mu }({U_{n + 1}}/{U_{n + 1}}I)\backslash {\mu }({U_n}/{U_n}I)$ there is $Q \in {\mu }({U_n}/{U_n}I)\backslash {\mu }({U_{n - 1}}/{U_{n - 1}}I)$ such that $Q < P$ and hence $d(Q) > d(P)$ and the relation (\ref{4.4}) follows. 
\par  
The relation (\ref{4.4}) easily implies that there is $m \in \mathbb{N}$ such that  
$$d({\mu }({U_{n + 1}}/{U_{n + 1}}I)\backslash {\mu }({U_n}/{U_n}I)) =  - 1$$ 
(i.e. ${\mu }({U_{n + 1}}/{U_{n + 1}}I)\backslash {\mu }({U_n}/{U_n}I) = \emptyset $) for all $n \geqslant m$. 
Therefore, $\mu(U_{n+1}/U_{n+1}I)\subseteq \mu (U_{n}/U_{n}I)$). Then it follows from ( \ref
{4.2}) that for any $t>n+1$ we have $\mu (U_{t}/U_{t}I)\cap \mu
(U_{n}/U_{n}I)=\mu (U_{n+1}/U_{n+1}I)$. Therefore, $\mu(U_{n+1}/U_{n+1}I)\subseteq \mu (U_{t}/U_{t}I)$. 
Then, as $t>n+1$ , it easily follows from (ii), that $\mu (U_{n+1}/U_{n+1}I)=\mu (U_{t}/U_{t}I)$  
but it contradicts (\ref{4.3}). 

\end{proof}


\begin{lemma}
Let $N$ be a minimax torsion-free nilpotent group and let $K$ be a normal subgroup of $N$ such that the quotient group $N/K$ is torsion-free abelian. 
Let $L$ be a dense subgroup of $N$ such that $K \leqslant L$ and the quotient group $L/K$ is free abelian.  Let $k$ be a field and let $W$ be a  $kN$-module which is  $kK$-torsion-free. Let $I$ be an $N$-large ideal  of $kK$ such that $W \ne WI$ and $A$ be a central subgroup of finite index in $L/{I^\dag }$. Then  there exists a cyclic $kN$-submodule $0 \ne U \leqslant W$ such that ${\mu _{kB}}(bkL/bkLI) = {\mu _{kB}}(dkL/dkLI)$ for any elements $0 \ne b,d \in U$ and any  subgroup $B\leq A$ of finite index in $A$.  
\end{lemma}

     \begin{proof} By Lemma 4.1(iv), it is sufficient to show that  there exists a cyclic $kN$-submodule $0 \ne U \leqslant W$ such that ${\mu _{kA}}(bkL/bkLI) = {\mu _{kA}}(dkL/dkLI)$ for any elements $0 \ne b,d \in U$.
\par      
      By Lemma 3.3(ii), there is a cyclic $kL$-submodule $0 \ne V \leqslant W$ such that for any  element $0 \ne a \in V$ and any element $0 \ne b \in akN$ the $kA$-module $bkL/bkLI$ has a finite series each of whose quotient is isomorphic to some section of the $kA$-module $akL/akLI$. Then it follows from \cite[Ch. II, \S 4, Proposition 16]{Bour} that 
\begin{equation}
          Sup{p_{kA}}bkL/bkLI \subseteq Sup{p_{kA}}akL/akLI.      \label{4.5}
\end{equation}
By Lemma 3.3(iii), there is a $kL$-submodule $0 \ne ckL \leqslant akL$ such that the $kA$-module $ckL/ckLI$ has a finite series each of whose quotient is isomorphic to some section of the $kA$-module $bkL/bkLI$. Then it follows from \cite[Ch. II, \S 4, Proposition 16]{Bour} that 
\begin{equation}
          Sup{p_{kA}}ckL/ckLI \subseteq Sup{p_{kA}}bkL/bkLI.      \label{4.6}
\end{equation}
By Lemma 4.2(iv), we also can assume that ${\mu _{kA}}({V_1}/{V_1}I) = {\mu _{kA}}(V/VI)$  for any cyclic $kL$-submodule  $0 \ne {V_1} \leqslant V$ and hence 
\begin{equation}
          {\mu _{kA}}(ckL/ckLI) = {\mu _{kA}}(akL/akLI)      \label{4.7}
\end{equation}
If $P \in {\mu _{kA}}(akL/akL)$ and $P \notin {\mu _{kA}}(bkL/bkLI)$ then it follows from (\ref{4.7}) that $P \in {\mu _{kA}}(ckL/ckL)$ and $P \notin {\mu _{kA}}(bkL/bkLI)$. Therefore, it follows from (\ref{4.6}) that there is $Q \in {\mu _{kA}}(bkL/bkLI)$ such that $Q < P$. But then (\ref{4.5}) implies that $P \notin {\mu _{kA}}(akL/akL)$ and the obtained contradiction shows that ${\mu _{kA}}(akL/akL) \subseteq {\mu _{kA}}(bkL/bkLI)$
\par  
If $P \in {\mu _{kA}}(bkL/bkLI)$ and $P \notin {\mu _{kA}}(akL/akL)$ then it follows from (\ref{4.5}) that there is $Q \in {\mu _{kA}}(akL/akLI)$ such that $Q < P$ and (\ref{4.7}) implies that $Q \in {\mu _{kA}}(ckL/ckLI)$. But then (\ref{4.6}) implies that $P \notin {\mu _{kA}}(bkL/bkL)$ and the obtained contradiction shows that ${\mu _{kA}}(akL/akL) \supseteq {\mu _{kA}}(bkL/bkLI)$. Thus, we can conclude that ${\mu _{kA}}(akL/akL) = {\mu _{kA}}(bkL/bkLI)$. Put $U=akN$ then we have 
${\mu _{kA}}(akL/akL) = {\mu _{kA}}(bkL/bkLI)$ for any element $0 \ne b \in U$ and the assertion follows. 
\end{proof}

	We also need some notions introduced by Wilson in \cite[Section 3.8]{Wils88} which are based on results of Brookes \cite{Broo85}. Let $A$ be a torsion-free abelian group of finite rank and let $k$ be a field. If $I$ and $J$ are ideals of $kA$ then we write $I \approx J$ if  $I \cap kB = J \cap kB$ for some finitely generated dense subgroup $B \leqslant A$. Then $ \approx $ is an equivalence relation on the set of all prime ideals of $kA$ and we denote by $\left[ I \right]$ the class of equivalence containing an ideal $I$. If a group $G$ acts on $A$ then we obtain an action of $G$ on the set of equivalent prime ideals of $kA$ which is given by ${\left[ I \right]^g} = \left[ {{I^g}} \right]$. 
\par  
	 If $B$ is a dense subgroup of $A$ and $P$ is a prime ideal of $kB$ then, as $kA$ is an integer domain over $kB$, it follows from \cite[Chap. V, \S 2, Theorem 1]{Bour} that there is a prime ideal $Q$ of $kA$ such  that $Q \cap kB = P$ and we put  ${[P]_{kA}} = [Q]$. If $\mu $ is a set of prime ideals of $kB$ then we put ${[\mu ]_{kA}} = \{ {[P]_{kA}}|P \in \mu \} $. 


\begin{lemma}
Let $A$ be a torsion-free abelian group of finite rank, let $B$ be a dense subgroup of $A$, let $k$ be a field and let  $I$ be an ideal of $kB$. Then: 
\begin{romanlist} 
\item there is a finitely generated dense subgroup $C \leqslant B$ such that for any subgroup $X$ of finite index in $C$ the mapping ${\mu _{kC}}(I \cap kC) \to {\mu _{kX}}(I \cap kX)$ given by $P \mapsto P \cap kX$  is  bijective;
\item for any finitely generated dense subgroup $C \leqslant B$ which meets (i) we have $|{[{\mu _{kC}}(I \cap kC)]_{kA}}| = |{\mu _{kC}}(I \cap kC)| < \infty $;   
\item  for any finitely generated dense subgroups $C,D \leqslant B$ which meet (i) we have ${[{\mu _{kC}}(I \cap kC)]_{kA}} = {[{\mu _{kD}}(I \cap kD)]_{kA}}$;
\end{romanlist}
\end{lemma}

\begin{proof} (i) It follows from \cite[Lemma 2.1]{Broo85} that there is a finitely generated dense subgroup $C \leqslant B$ such that for any subgroup $X$ of finite index in $C$ if $P \in {\mu _{kX}}(I \cap kC)$ then ${\varphi _X}:P \mapsto P \cap kX$ maps ${\mu _{kC}}(I \cap kC)$ into ${\mu _{kX}}(I \cap kX)$. It follows from Lemma 4.1(i) that the mapping ${\varphi _X}$ is surjective. If ${\varphi _X}$ is not injective for some subgroup $X$ of finite index in $C$, then $|{\mu _{kX}}(I \cap kX)| < |{\mu _{kC}}(I \cap kC)|$ and we can replace the subgroup $C$ by its subgroup $X$ of finite index with minimal possible $|{\mu _{kX}}(I \cap kX)|$. 
\par  
(ii) Let $[P],[Q] \in {[{\mu _{kC}}(I \cap kC)]_{kA}}$ then $P \cap kC$ and $Q \cap kC$ are in ${\mu_{kC}}(I \cap kC)$. $[P] = [Q]$ if and only if there is a finitely generated dense subgroup $X$ of $A$ such that $P \cap kX = Q \cap kX$. Evidently, we can assume that $X \leqslant C$ and, as  the mapping ${\mu _{kC}}(I \cap kC) \to {\mu _{kX}}(I \cap kX)$ given by $P \mapsto P \cap kX$  is  bijective, we can conclude that $P \cap kC = Q \cap kC$. Therefore, 
for any $U,V \in {\mu _{kC}}(I \cap kC)$ if $U \ne V$ then ${[U]_{kA}} \ne {[V]_{kA}}$ and it implies that 
$\left| {{{[{\mu _{kC}}(I \cap kC)]}_{kA}}} \right| = \left| {{\mu _{kC}}(I \cap kC)} \right|$. 
As the ring $kC$ is Noetherian, it follows from \cite[Chap. IV, \S 1, Theorem 2]{Bour} that $|{\mu _{kC}}(I \cap kC)| < \infty $.
\par  
(iii) It is easy to note that the subgroup $X = D \cap C \leq A$ meets (i) and hence it is sufficient to show that ${[{\mu _{kC}}(I \cap kC)]_{kA}} = {[{\mu _{kX}}(I \cap kX)]_{kA}} = {[{\mu _{kD}}(I \cap kD)]_{kA}}$. Thus, we can consider only the case, where $X = D \leqslant C$. If $[P] \in {[{\mu _{kC}}(I \cap kC)]_{kA}}$ then, by the definition of ${[{\mu _{kC}}(I \cap kC)]_{kA}}$, we can assume that $P \cap kC \in {\mu _{kC}}(I \cap kC)$. As by (i), ${\varphi _X}:P \mapsto P \cap kX$ maps ${\mu _{kC}}(I \cap kC)$ into ${\mu _{kX}}(I \cap kX)$ we have $P \cap kX \in {\mu _{kX}}(I \cap kX)$ and hence $[P] \in {[{\mu _{kX}}(I \cap kX)]_{kA}}$. On the other hand, if $[P] \in {[{\mu _{kX}}(I \cap kX)]_{kA}}$ then, by the definition of ${[{\mu _{kX}}(I \cap kX)]_{kA}}$, we can assume that $P \cap kX \in {\mu _{kX}}(I \cap kX)$. By Lemma 4.1(i), there is $Q\in {\mu _{kC}}(I \cap kC)$ such that $Q\cap kX = P\cap kX$. Let $[S]=[Q]_{kA}$ then $[S]\in [\mu_{kC}(I \cap kC)]_{kA}$  and $S\cap kX = P\cap kX$. Therefore, $[P]= [S]\in [\mu_{kC}(I \cap kC)]_{kA}$.
\end{proof}
 
Let $G$ be a group and let $K$ be a normal subgroup of $G$. Let $N$ be a nilpotent subgroup of $G$ such that $K \leqslant N \leqslant G$ and the quotient group $N/K$ is torsion-free abelian of finite rank. Let $k$ be a field and let $I$ be a $G$-large ideal of $kK$. Then, as the quotient group $K/I^\dag$ is finite and the quotient group $N/K$ is torsion-free abelian of finite rank, it follows from Lemma 2.2(i) that 
$N/I^\dag$ has a characteristic central torsion-free subgroup $A$ of finite index. 
\par  
Let $L$ be a dense subgroups of $N$ such that $K \leqslant L$ and the quotient groups $L/K$ is finitely generated. Then $A \cap L/I^\dag$ 
is a dense central finitely generated torsion-free subgroup of $A$. 
\par  
Let $W$ be a $kN$-module which is $kK$-torsion-free and let $akL$ be a cyclic $kL$-module generated by an element $0 \ne a \in W$. 
By Lemma 4.4(i), there is a finitely generated dense subgroup $A_L \leqslant A \cap L/I^\dag$ 
such that for any subgroup $X$ of finite index in ${A_L}$ the mapping ${\mu _{k{A_L}}}(An{n_{k{A_L}}}(akL/akLI)) \to {\mu _{kX}}(An{n_{kX}}(akL/akLI))$  given by $P \mapsto P \cap kX$ is bijective. 
\par  
As ${A_L}$ is a finitely generated  subgroup of finite index in $A \cap L / I^\dag$, we can conclude that ${A_L}$ is a dense finitely generated  subgroup of $A$ and $akL/akLI$ is a finitely generated $k{A_L}$-module. Then it follows from \cite[Theorem 3.7]{Wehr09} that the domain $k{A_L}$ is Noetherian and $akL/akLI$ is a Noetherian $k{A_L}$-module. Thus, the sets $Sup{p_{k{A_L}}}(akL/akLI)$ and ${\mu _{k{A_L}}}(akL/akLI)$ are well defined. Then according to (\ref{4.1}) we have  
\begin{equation}
                         {\mu _{k{A_L}}}(akL/akLI) = {\mu _{k{A_L}}}(An{n_{k{A_L}}} (akL/akLI))                 \label{4.8} 
\end{equation}
Thus, the set ${\mu _{k{A_L}}}(akL/akLI)$ is defined for any $0 \ne a \in W$ and we put ${[{\mu _{k{A_L}}}(akL/akLI)]_{kA}} = \left\{ {{{[P]}_{kA}}|P \in {\mu _{k{A_L}}}(akL/akLI)} \right\}$. It follows from (\ref{4.8}) that 
\begin{equation}
             {[{\mu _{k{A_L}}}(akL/akLI)]_{kA}} = {[{\mu _{k{A_L}}}(An{n_{k{A_L}}}(akL/akLI))]_{kA}}           \label{4.9}   
\end{equation}
By (\ref{4.9}) and according to Lemma 4.4(ii),(iii), the set ${[{\mu _{k{A_L}}}(akL/akLI)]_{kA}} = \left\{ {{{[P]}_{kA}}|P \in {\mu _{k{A_L}}}(An{n_{k{A_L}}}({a_L}))} \right\}$ is finite and does not depend on the choice of the subgroup ${A_L}$ which meets the conditions of Lemma 4.4(i). Everywhere below in the definition of the set 
$${[{\mu _{k{A_L}}}(akL/akLI)]_{kA}} = \left\{ {{{[P]}_{kA}}|P \in {\mu _{k{A_L}}}(akL/akLI)} \right\}$$ 
we assume that the subgroup ${A_L}$ meets the conditions of Lemma 4.4(i),i.e. for any subgroup $X$ of 
finite index in ${A_L}$ the mapping ${\mu _{k{A_L}}}(An{n_{k{A_L}}}(akL/akLI)) \to {\mu _{kX}}(An{n_{kX}}(akL/akLI))$ 
 given by $P \mapsto P \cap kX$ is bijective. .


\begin{proposition}
Let $N$ be a minimax torsion-free nilpotent group and let $K$ be a normal subgroup of $N$ such that the quotient group $N/K$ is torsion-free abelian. Let $k$ be a field such that $chark = 0$ and let $W$be a $kN$-module which is  $kK$-torsion-free. Let $I$ be an $N$-large ideal of $kK$ and $A$ be a torsion-free characteristic central subgroup of finite index in $N / I^\dag$. Then there exists a cyclic $kN$-submodule $0 \ne V \leqslant W$ such that ${[{\mu _{k{A_L}}}(akL/akLI)]_{kA}} = {[{\mu _{k{A_M}}}(bkM/bkMI)]_{kA}}$ for any elements $0 \ne a,b \in V$ and any dense subgroups $L,M \leqslant N$ such that $K \leqslant L \cap M$ and the quotient groups $L/K$ and $M/K$ are finitely generated. 
\end{proposition}

\begin{proof} At first we fix a dense subgroup $L \leqslant N$ such that  $K \leqslant L$ and the quotient group $L/K$ is finitely generated. By Lemma 4.3, we can choose an element $0 \ne a \in W$ such that ${\mu _{k{B_L}}}(akL/akLI) = {\mu _{k{B_L}}}(bkL/bkLI)$ for any $0 \ne b \in akN$  and any subgroup $B_L$ of finite  index in $A_L$. Put $V = akN$ then it is sufficient to show that ${[{\mu _{k{A_L}}}(bkL/bkLI)]_{kA}} = {[{\mu _{k{A_M}}}(bkM/bkMI)]_{kA}}$ for any $0 \ne b \in akN$ and any dense subgroup $M \leqslant N$ such that $K \leqslant M$ and the quotient group $M/K$ is finitely generated. 
\par  
It is easy to note that $T = L \cap M$ is a subgroup of finite index in $L$ and $M$ . It easily implies that  we can choose a finitely generated dense subgroup  ${A_T} \leqslant A \cap (T/{I^\dag })$ such that${A_T} \leqslant {A_L} \cap {A_M}$. Then it is sufficient to show that ${[{\mu _{k{A_L}}}(bkL/bkLI)]_{kA}} = {[{\mu _{k{A_T}}}(bkT/bkTI)]_{kA}} = {[{\mu _{k{A_M}}}(bkM/bkMI)]_{kA}}$. In fact it is sufficient to prove the identity $${[{\mu _{k{A_L}}}(bkL/bkLI)]_{kA}} = {[{\mu _{k{A_T}}}(bkT/bkTI)]_{kA}}$$ 
because the identity ${[{\mu _{k{A_T}}}(bkT/bkTI)]_{kA}} = {[{\mu _{k{A_M}}}(bkM/bkMI)]_{kA}}$ is analogous. 
\par  
By Lemma 4.4(iii), ${[{\mu _{k{A_L}}}(bkL/bkLI)]_{kA}} = {[{\mu _{k{A_T}}}(bkL/bkLI)]_{kA}}$ and hence it is sufficient to show that ${[{\mu _{k{A_T}}}(bkT/bkTI)]_{kA}} = {[{\mu _{k{A_T}}}(bkL/bkLI)]_{kA}}$. Since $T \leqslant L$, we have $bkT \leqslant bkL$ and it  follows from Lemma 4.2(i) that 
\begin{equation}
          Sup{p_{k{A_T}}}bkT/bkTI \subseteq Sup{p_{k{A_T}}}bkL/bkLI.                      \label{4.10}
\end{equation}
On the other hand, as $T \leqslant L$, we see that $bkL = \sum\limits_{i = 1}^n {bkT{g_i}} $ for some ${g_i} \in L$. Therefore, 
\begin{equation}
                    bkL/bkLI \cong (\sum\limits_{i = 1}^n {bkT{g_i}} /bkT{g_i}I)/X,                  \label{4.11} 
\end{equation}
where $X$ is a $k{A_T}$-submodule of $\sum\limits_{i = 1}^n {bkT{g_i}} /bkT{g_i}I$. As ${A_T}$ is a central subgroup of $N/{I^\dag }$, we can conclude that $bkT{g_i}/bkT{g_i}I \cong bkT/bkTI$ for all ${g_i} \in L$. Then it follows from (\ref{4.11}) and \cite[Ch. II, \S 4, Proposition 16]{Bour} that 
\begin{equation}
                 Sup{p_{k{A_T}}}bkL/bkLI \subseteq Sup{p_{k{A_T}}}bkT/bkTI.                 \label{4.12} 
\end{equation}
The relations (\ref{4.10}) and (\ref{4.12}) imply that $Sup{p_{k{A_T}}}bkL/bkLI = Sup{p_{k{A_T}}}bkT/bkTI$ and hence ${\mu _{k{A_T}}}(bkL/bkLI) = {\mu _{k{A_T}}}(bkT/bkTI)$.  Therefore, ${[{\mu _{k{A_T}}}(bkT/bkTI)]_{kA}} = {[{\mu _{k{A_T}}}(bkL/bkLI)]_{kA}}$. 
\end{proof}


\begin{corollary}
In the denotations of Proposition 4.1, the cyclic $kN$-submodule $0 \ne V \leqslant W$ defines a finite set
\begin{equation}
 {M_{kA}}(V/VI) = {[{\mu _{k{A_L}}}(akL/akLI)]_{kA}}  \label{4.13}
 \end{equation}
 of equivalent classes of prime ideals of $kA$ which depends only on the ideal $I$ and the subgroup $A$. 
\end{corollary}

\begin{proof} The finiteness on (\ref{4.13}) follows from (\ref{4.9})  and Lemma 4.4(ii),(iii). By Proposition 4.1, (\ref{4.13}) does not depend on the choice of the subgroup $L$ and the element $0 \ne a \in V$. By Lemma 4.4(iii), 
(\ref{4.13}) does not depend on the choice of the subgroup ${A_L}$ which meets the conditions of Lemma 4.4(i). 
\end{proof}

Combining Propositions 3.1, 3.2 and Corollary 4.1 we obtain the following theorem


\begin{theorem}
Let $G$ be a soluble group of finite torsion-free rank ${r_0}(G) < \infty $, let $N$ be a minimax nilpotent torsion-free normal subgroup of $G$ and let $K$ be a $G$-invariant subgroup of $N$ such that the quotient group $N/K$ is torsion-free abelian. Let $L$ be a dense subgroup of $N$ such that $K \leqslant L$ and the quotient group $L/K$ is free abelian and let $\chi $ be a full system of subgroups of $L$ over $K$. Let $k$ be a field of characteristic zero and let $W$ be a $kN$-module which is $kK$-torsion-free. Then there are a cyclic $kN$-submodule $0 \ne V \leqslant W$, a $G$-large ideal $I$ of $kK$ and a central $G$-invariant subgroup $A$ of finite index in $N/{I^\dag }$ such that for any element $0 \ne b \in V$ we have:     
\begin{romanlist}
\item $X \in \chi (bkL)$ if and only if  the quotient module $bkL/bkLI$ is not $k{A_X}$-torsion, where ${A_X} = A \cap X/{I^\dag }$;    
\item for any dense subgroup $M \leqslant N$ with $K \leqslant M$ and finitely generated quotient group $M/K$ the finite set ${M_{kA}}(V/VI) = {[{\mu _{k{A_M}}}(bkM/bkMI)]_{kA}}$ of equivalent classes of prime ideals of $kA$ depends only on the ideal $I$ and the subgroup $A$;
 \item for any  $g \in G$ we have ${M_{kA}}(Vg/VgI) = {  M_{kA}}{(V/VI)^g} = \{ {[P]^g} = [{P^g}]|[P] \in {M_{kA}}(V/VI)\} $. 
\end{romanlist}
\end{theorem}

	\begin{proof}  (i) By Propositions 3.2, there are a cyclic $kN$-submodule $0 \ne U \leqslant W$, a finitely generated dense subgroup $H$ of $L$ and a right ideal $J$ of $k(H \cap K)$  such that if $I$ is an $N$-large ideal of $kK$ such that $I \cap k(K \cap H)$ culls $J$ in $k(H \cap K)$ then  for any $0 \ne b \in U$ we have $X \in \chi (bkL)$ if and only if  the quotient module $bkL/bkLI$ is not $k{A_X}$-torsion, where ${A_X} = A \cap X/{I^\dag }$ and $A$  is central subgroup of finite index in $N/{I^\dag }$. By Propositions 3.1, there is a $G$-large ideal $I$ of $kK$ such that $I \cap k(K \cap H)$ culls $J$ in $k(H \cap K)$. As $N \leqslant G$, we can conclude that the ideal $I$ is an $N$-large and hence the assertion (i) holds for the ideal $I$ and the submodule $U$. Since $N$, ${I^\dag }$ and $K$ are normal subgroups of $G$ and any abelian-by-finite nilpotent group of finite rank has a characteristic central subgroup of finite index, we can chose the subgroup $A$ such that $A$ is $G$-invariant. Then replacing $W$ by $U$ we can assume that (i) holds for any $0 \ne b \in W$. 
\par  
(ii) By Corollary 4.1, we can choose a cyclic submodule $0 \ne V \leqslant W$ such that for any $0 \ne b \in V$ and any dense subgroup $M \leqslant N$ with $K \leqslant M$ and finitely generated quotient group $M/K$ the finite set ${ M_{kA}}(V/VI) = {[{\mu _{k{A_M}}}(bkM/bkMI)]_{kA}}$ depends only on the ideal $I$ and the subgroup $A$. 
\par  
(iii) Since the subgroup $A$ and the ideal $I$ are $G$-invariant, we have $ M_{kA}(Vg/VgI) = [{\mu _{kA_L^g}}(bgkL/bgkLI)]_{kA}$. 
Then by (\ref{4.9}), $M_{kA}(Vg/VgI)=[\mu _{kA_L^g}(bgkL/bgkLI)]_{kA} = {[{\mu _{kA_L^g}}(An{n_{kA_L^g}}(bgkL/bgkLI))]_{kA}}$. It is not difficult to show that  
$$[\mu_{kA_L^g}(Ann_{kA_L^g}(bgkL/bgkLI))]_{kA}={[{\mu _{kA_L}}({(An{n_{k{A_L}}}(bkL/bkLI))^g})]_{kA}}$$  
and hence $M_{kA}(Vg/VgI)={[{\mu _{kA_L}}({(An{n_{k{A_L}}}(bkL/bkLI))^g})]_{kA}}$.

Thus, $M_{kA}(Vg/VgI)=[{({\mu _{k{A_L}}}{(Ann_{kA_L}(bkL/bkLI))^g}]_{kA}}=\{ {[{P^g}]_{kA}}|P \in ({\mu _{k{A_L}}}(An{n_{k{A_L}}}(bkL/bkLI))\} = \{ {[P]^g}_{kA}|P \in ({\mu _{k{A_L}}}(An{n_{k{A_L}}}(bkL/bkLI))\}={\{ [P]^g = [P^g ]|[P] \in M_{kA}(V/VI)\}= M_{kA}}{(V/VI)^g}$. 
\end{proof}


\section{Uniform Stabilized Modules over Minimax Nilpotent Groups}

      Let $A$  be an abelian torsion-free group of finite rank acted by a group $G$ , let $k$ be a field and let $I$ be an ideal of $kA$. A subgroup ${S_G}(I) \leqslant G$ which consists of all elements  $g \in G$ such that $I \cap kB = I^g  \cap kB$ for some finitely generated dense subgroup $B \leqslant A$ is said to be the standardizer of $I$ in $G$ (see \cite{Broo85,Wils88}). Since the class $[I]$ consists of all ideals $J \leqslant kA$ such that that $I \cap kB = J \cap kB$ for some finitely generated dense subgroup $B \leqslant A$, we see that $[I]^g = \left\{ J^g |J \in [I] \right\}$ also forms a class $[{I^g}]$. So, we have an action of $G$ on the set of classes of ideals of $kA$. It immediately follows from the definition of ${S_G}(I)$ that ${S_G}(I) = {S_G}([I])$, where ${S_G}([I]) = \left\{ {g \in G|{[I]^g} = [I]} \right\}$ is the stabilizer of $[I]$ in $G$. 

\begin{lemma}
Let $R$ be a Noetherian domain, let $A$ be a finitely generated abelian group and let $I$ be an ideal of the group ring $RA$. Let $B$ be a torsion-free subgroup of $A$ such that $RB \cap P \ne 0$ for any $P \in {\mu _R}(I)$. Then $RB \cap I \ne 0$. 
\end{lemma}
	\begin{proof} By \cite[Theorem 3.7]{Wehr09}, the ring $RA$ is Noetherian and it follows from \cite[Corollary 37.11]{Pass89} that $kB$ is a domain. Then, by \cite[Chap. IV, \S 5, Theorem 5]{ZaSa58}, the ideal $I$ has a presentation $I =  \cap _{i = 1}^n{D_i}$, where ${D_i}$ are primary ideals of $RA$. Let ${P_i}$ be the prime radical of ${D_i}$, then ${P_i}$ is a prime over $I$ and it follows from the definition of ${\mu _R}(I)$ that there is $P \in {\mu _R}(I)$ such that ${P_i} \geqslant P$. So, as $RB \cap P \ne 0$ , we can conclude that there is an element $0 \ne {a_i} \in RB \cap {P_i} \ne 0$ for each $i$. Then it follows from \cite[Chap. III, \S 9, Theorem 13]{ZaSa58} that that there is an integer ${m_i} \in \mathbb{N}$ such that ${a_i}^{{m_i}} \in {D_i}$ and therefore $0 \ne \Pi _{i = 1}^n{a_i}^{{m_i}} \in I \cap RB$.
\end{proof}


\begin{lemma}
Let $A$ be an abelian torsion-free group of finite rank acted by a soluble-by-finite group $G$, let $C = So{c_G}A$ and let $B$ be a finitely generated dense subgroup of $A$. Let $k$ be a field of characteristic zero and let $J$ be an ideal of $kB$. Let ${P_1},{P_2},...,{P_n}$ be prime ideals of $kA$ such that ${\mu _{kB}}(J) = \left\{ {{P_i} \cap kB|1 \leqslant i \leqslant n} \right\}$ is the set of minimal prime ideals over $J$. If $\left| {G:{S_G}([{P_i}])} \right| < \infty $ for all $1 \leqslant i \leqslant n$ then $J \cap k(C \cap B) \ne 0$. 
\end{lemma}
\begin{proof} Since $\left| {G:{S_G}([{P_i}])} \right| < \infty $ for all $1 \leqslant i \leqslant n$, we can conclude that $G$ has a normal subgroup $H$ of finite index such that $H = {S_H}([{P_i}])$ for all $1 \leqslant i \leqslant n$. By Lemma 2.4(i), $Soc_{H}A = Soc_{G}A$ and hence there is no harm in assuming that $H = G$. Thus, we may assume that $G = {S_G}([{P_i}]) = {S_G}({P_i})$ for all $1 \leqslant i \leqslant n$. 
\par  
	Suppose that ${P_i}^\dag  = \left\{ {a \in A|1 - a \in {P_i}} \right\} \ne 1$. As $is{({P_i}^\dag )^g} = is({({P_i}^g)^\dag })$ for all $g \in G$ and $G = {S_G}({P_i})$, we can conclude that $is{({P_i}^\dag )^g} = is({P_i}^\dag )$ for all $g \in G$, i.e. $is({P_i}^\dag )$ is a $G$-invariant subgroup of $A$. Then $is({P_i}^\dag ) \cap C \ne 1$ and as $C$ is an isolated subgroup of $A$, we see that ${P_i}^\dag  \cap C \ne 1$. It easily implies that ${P_i} \cap kC \ne 0$. 
\par  
	Suppose now that ${P_i}^\dag  = \left\{ {a \in A|1 - a \in {P_i}} \right\} = 1$. Then by \cite[Theorem A]{Broo85} (see also \cite[Theorem 2.4]{Tush14}), ${P_i} = ({P_i} \cap k{\Delta_G}(A))kA$ and hence ${P_i} \cap kC \ne 0$ because, by Lemma 2.3(i), $\Delta_{G}(A) \leqslant Soc_{G}(A) = C$. 
\par  
	Thus, ${P_i} \cap kC \ne 0$ for all $1 \leqslant i \leqslant n$. Then, as $B$ is a dense subgroup of $A$, it follows from \cite[Lemma 2.2.3(ii)]{Tush2000} that ${P_i} \cap k(C \cap B) \ne 0$ for all $1 \leqslant i \leqslant n$. Since ${\mu _{kB}}(J) = \left\{ {{P_i} \cap kB|1 \leqslant i \leqslant n} \right\}$ and ${P_i} \cap k(C \cap B) \ne 0$ for all $1 \leqslant i \leqslant n$, it follows from Lemma 5.1 that $J \cap k(C \cap B) \ne 0$. 
\end{proof}

Let $R$ be a ring and let $M$, $X$ and $Y$ be $B$-modules. The modules $X$ and $Y$ are separated in  $M$ if $X$ and $Y$ do not have non-zero isomorphic sections which are isomorphic to a submodule of $M$. Let $k$ be a field , $G$ be a group and $A$ be a normal subgroup of $G$. Let $M$ and $W$ be $kA$-modules. The subgroup $Sep_{(G,A)}(M,W) \leqslant G$ generated by all elements $g \in G$ such that $kA$-modules $W$ and $Wg$ are not separated in $M$ is said to be the separator of $W$ in $G$. Evidently, for any element $h \in G$ such that $h \notin Sep_{(G,A)}(M,W)$ modules $W$ and $Wh$ - are separated in $M$. 
\par  
$kA$-modules $W$ and $V$ are said to be similar if their injective hulls $[W]$ and $[V]$ are isomorphic. The modules $W$ and $V$ are similar if and only if they have isomorphic essential submodules. By \cite[Lemma 3.2]{BBro85}, the stabilizer $Stab_G W = \{ g \in G | Wg \ and \ W \ are \ similar \}$ of $W$ in $G$ is a subgroup of $G$. 
\par  
We say that a submodule $X \leqslant M$ is solid in $M$ if $X$ is uniform and $M$ do not have submodules which are isomorphic to a proper section of $X$. It is easy to note that $Stab_G W \leqslant Sep_{(G,A)}(M,W)$, moreover, if the submodule $W$ is solid in $M$ then $Stab_G W = Sep_{(G,A)}(M,W)$. 


\begin{proposition}
Let $N$ be a nilpotent normal non-abelian minimax torsion-free subgroup of a solvable-by-finite group $G$ such that ${r_0}(G) < \infty $. Let $K$ be an isolated $G$-invariant subgroup of $N$ such that the quotient group $N/K$ is abelian. Let $k$ be a field such that $chark = 0$. Suppose that there is an uniform $kN$-torsion $kN$-module $W$ such that $Stab_G W = G$ and for any proper isolated $G$–invariant subgroup $X$ of $N$ such that $K \leqslant X$ the module $W$ is $kX$-torsion-free. Then $Soc_G N/K = N/K$. 
\end{proposition}

\begin{proof} Put $Soc_{G} N/K = S/K$. Suppose that $Soc_{G}N/K \ne N/K$ then it follows from Lemma 2.3(i) that $S$ is a proper isolated $G$-invariant subgroup of $N$ such that $K \leqslant S$ and hence the module $W$ is $kS$-torsion-free. 
\par  
Let $L$ be a dense subgroup of $N$ such that $K \leqslant L$ and the quotient group $L/K$ is free abelian. As the subgroup $S$ is isolated in $N$, we see that $T = L \cap S$ is isolated in $L$ and hence we can choose a full system 
$\chi  = \{ \left\langle K,\{ x_j|j \in J \} \right\rangle |J \subseteq \{ 1,...,n \} \} $ 
of subgroups of $L$ over K such that $T \in \chi $. As the module $W$ is $kS$-torsion-free, the module $W$ is also $kT$-torsion-free. 
\par  
By Theorem 4.1, there are a cyclic $kN$-submodule $0 \ne V \leqslant W$, a $G$-large ideal $I$ of $kK$ and a central $G$-invariant subgroup $A$ of finite index in $N/{I^\dag }$ such   that for any element $0 \ne b \in V$ we have: 
\begin{romanlist}
\item  $X \in \chi (bkL)$ if and only if  the quotient module $bkL/bkLI$ is not $k{A_X}$-torsion, where ${A_X} = A \cap X/{I^\dag }$ and $A$  is central subgroup of finite index in $N/{I^\dag }$;    
\item for any dense subgroup $M \leqslant N$ with $K \leqslant M$ and finitely generated quotient group $M/K$ the finite set ${M_{kA}}(V/VI) = {[{\mu _{k{A_M}}}(bkM/bkMI)]_{kA}}$ depends only on the ideal $I$ and the subgroup $A$;
\item for any  $g \in G$ we have ${M_{kA}}(Vg/VgI) = {M_{kA}}{(V/VI)^g} = \{ {[P]^g} = [{P^g}]|[P] \in {M_{kA}}(V/VI)\} $. 
\end{romanlist}

	Since the module $W$ is $kT$-torsion-free, we have $T \in \chi (bkL)$ for any element $0 \ne b \in V$. As the module $W$ is uniform, we can conclude that $[V] = [W]$ and $[Vg] = [Wg]$ for any $g \in G$ and hence $[Vg] \cong [V]$. It implies that $Sta{b_G}V = G$. Then for any $g \in G$ there are non-zero elements $b,c \in V$ such that $bkN \cong ckNg = cgkN$. Put $M = L \cap {L^{ - g}}$ then it follows from (ii) that ${M_{kA}}(V/VI) = {[{\mu _{k{A_M}}}(bkM/bkMI)]_{kA}} = {[{\mu _{k{A_M}}}(bgkM/bgkMI)]_{kA}} = $ ${M_{kA}}(Vg/VgI)$. Therefore, as by (iii) ${M_{kA}}(Vg/VgI) = {M_{kA}}{(V/VI)^g}$, we see that ${M_{kA}}(V/VI) = {M_{kA}}{(V/VI)^g}$ = $\left\{ {{{[P]}^g} = [{P^g}]|[P] \in {M_{kA}}(V/VI)} \right\}$ and hence the group $G$ acts on the finite set ${M_{kA}}(V/VI)$. Then $\left| {G:{S_G}([P])} \right| < \infty $ for any $[P] \in {M_{kA}}(V/VI)$ and it easily implies that there is a subgroup $H$ of finite index in $G$ such that ${[P]^g} = [{P^g}] = [P]$ for any $[P] \in {M_{kA}}(V/VI)$ and any $g \in H$.  
\par  
By Corollary 4.1, we have ${M_{kA}}(V/VI) = {[{\mu _{k{A_L}}}(bkL/bkLI)]_{kA}}$. Suppose that ${M_{kA}}(V/VI) = \left\{ {[{P_1}],...,[{P_n}]} \right\}$ then $B = {A_L}$ is a finitely generated dense subgroup of $A$ such that ${\mu _{kB}}(An{n_{kB}}(bkL/bkLI)) = \left\{ {{P_i} \cap kB|1 \leqslant i \leqslant n} \right\}$. Put $C = A \cap S/{I^\dag }$, as $So{c_G}N/K = S/K$, we see that $C = So{c_G}A$. Then, by Lemma 5.2, $An{n_{kB}}(bkL/bkLI) \cap k(C \cap B) \ne 0$ and hence $An{n_{kB}}(bkL/bkLI) \cap kC \ne 0$. If ${A_T} = A \cap T/{I^\dag }$, where $T = L \cap S$, then ${A_T}$ is a finitely generated dense subgroup of $C = A \cap S/{I^\dag }$ and hence, as $An{n_{kB}}(bkL/bkLI) \cap kC \ne 0$, it follows from \cite[Lemma 2.2.3(ii)]{Tush2000} that $An{n_{kB}}(bkL/bkLI) \cap k{A_T} \ne 0$. It implies that the module $bkL/bkLI$ is $k{A_T}$-torsion. But it contradicts (i) because $T \in \chi (bkL)$. Thus, the obtained contradiction shows that $So{c_G}N/K = N/K$. 
\end{proof}
     Let $R$ be a ring and let $M$ be an $R$-module then ${K_R}(M)$ denotes the Krull dimension of $M$. The module $M$ is said to be $\rho $-critical  if ${K_R}(M) = \rho $ and ${K_R}(M/U) < {K_R}(M) = \rho $ for any non-zero submodule $U \leqslant M$.
\par  
Let $N$ be a nilpotent group such that ${r_0}(N) < \infty $, let $k$ be a field and let $M$ be an $kN$-module.  Let $S$ be a finitely generated subring of $k$. Let $0 \ne a \in M$ and let $H$ be a proper isolated normal subgroup of $N$. We say that $(a,H)$ is an important pair for the $SN$-module $M$ if there is a finitely generated dense subgroup $A \leqslant N$ such that:
\begin{romanlist}
\item the module $aSA$ is $\rho $-critical  and ${K_{SX}}(aSX)$$ \leqslant $${K_{SX}}(xSX)$ for any element 0$ \ne $$x$$ \in $$M$ and any finitely generated dense subgroup  $X$$ \leqslant $$N$; 
 \item $aSA$=$aSB{ \otimes _{SB}}SA$, where $B$=$A$$ \cap $$H$; 
\item if 0$ \ne $$b$$ \in $$aSB$ and $V$ is a dense subgroup of $B$ then there is no isolated subgroup $D$$ \leqslant $$V$ such that $bSV = bSD{ \otimes _{SD}}SV$ and $i{s_N}(D)$ is a normal subgroup of $N$. 
\end{romanlist}
\par  
	The module $M$ is said to be impervious if it has no important pairs for any finitely generated subring $S$ of  $k$. 


\begin{proposition}
Let $N$ be a nilpotent  non-abelian minimax torsion-free normal  subgroup of a soluble group $G$ such that ${r_0}(G) < \infty $. Let $K$ and $L$ be isolated  $G$-invariant subgroups of $N$ such that $L \leqslant K$, $K$ contains next to the last member of the upper central series of $N$ and the quotient group $N/L$ is abelian. Let $k$ be a field such that $chark = 0$. Let $W$ be a $kG$-module  which is $kN$-torsion and for any proper isolated $G$–invariant subgroup $X$ of $N$ such that $L \leqslant X$ the module $W$ is $kX$-torsion-free. Suppose that $W$ is impervious $kN$-module. If $Se{p_{(G,Y)}}(xkG,xkY) = G$ for any element $0 \ne x \in W$ and any $G$-invariant subgroup $Y$ of $N$ then:
\begin{romanlist}
\item $So{c_G}N/L = N/L$; 
\item there is a dense $G$-invariant subgroup $D \leqslant N$ such that $K \leqslant D$ and the quotient group $D/K$ is polycyclic.     
\end{romanlist}
\end{proposition}

\begin{proof} Evidently, there is no harm in assuming that the module $W = akG$ is cyclic. Then, as the module $W = akG$ is $kN$-torsion, there is a finitely generated subfield $f \leqslant k$ such $J = An{n_{fN}}(a) \ne 0$ and hence, by Lemma 2.5, the module $W = akG$ is $fN$-torsion. Since for any proper isolated $G$–subgroup $X$ of $N$ such that $K \leqslant X$ the module $W$ is $kX$-torsion-free, it is also $fX$-torsion-free. By \cite[Lemma 3.1.2]{Tush2000}, $Se{p_{(G,kN)}}(W,xkN) \leqslant Se{p_{(G,fN)}}(W,xfN)$ and, as $Se{p_{(G,N)}}(W,xkN) = G$, we may conclude that $Se{p_{(G,N)}}(W,xfN) = G$. As the $kN$-module $xkN$ is impervious, it is easily to note that the $fN$-module $xfN$ is also impervious. Thus, there is no harm in assuming that the field $k$ is finitely generated. Then, by \cite[Lemma 3.1.2]{Tush2000}, we can choose the element $x$ such that the $kN$-module $xkN$ is solid in $W$and hence, the module $xkN$ is uniform and $Sta{b_G}xkN = Se{p_{(G,N)}}(W,xkN) = G$. 
\par  
          (i) Since the module $xkN$ is uniform and $Sta{b_G}xkN = Se{p_{(G,N)}}(W,xkN) = G$ , by Proposition 5.1, we have $So{c_G}N/K = N/K$. \par

(ii) It follows from (i) and Lemma 2.4(ii) that $G$ has a normal subgroup $H$ of finite index such that $So{c_G}N/K = N/K$ contains a dense $G$-invariant Noetherian $H$-polyplinth $Y/K$. As the subgroup $Y$ is dense in $N$ and the module $W$ is $kN$-torsion, it follows from \cite[Lemma 2.2.3(ii)]{Tush2000} that the module $W$ is $kY$-torsion. Since $W$ is impervious as a $kN$-module, $W$ is impervious as a $kY$-module. Then, as $Se{p_{(G,Y)}}(xkG,xkY) = G$, there is no harm in changing $N$ by $Y$, and we can assume that $A = N/K$ is a Noetherian $H$-polyplinth for some normal subgroup $H$ of finite index in $G$. After these changings, by \cite[Lemma 3.1.2]{Tush2000}, we can choose the element $x$ such that the $kN$-module $xkN$ is solid in $W$ and hence, the module $xkN$ is uniform and $Sta{b_G}xkN = Se{p_{(G,N)}}(W,xkN) = G$. By Lemma 2.4(iii), ${\Delta _G}(A) = D/K$ is a finitely generated isolated $G$-invariant subgroup of $A = N/K$ and, by Lemma 2.3(i), $\Delta_{G}(A/\Delta_{G}(A))$ is trivial. Then, as the group $D/K$ is polycyclic, changing $K$ by $D$ we can assume that ${\Delta _G}(A)$ is trivial. Since $Stab_{G}xkN = Sep_{(G,N)}(W,xkN) = G$, we have $Stab_{H}xkN = H$. Then it follows from \cite[Lemma 3.3.4]{Tush2000}  that the quotient group $N/K$ is polycyclic. 
\end{proof}


\begin{proposition}
Let $N$ be a nilpotent non-abelian minimax torsion-free normal  subgroup of a soluble group $G$ with ${r_0}(G) < \infty $. Let $k$ be a field such that $chark = 0$.  Let $W$ be a $kG$-module  which is $kN$-torsion and for any proper isolated $G$–invariant subgroup $X$ of $N$ the module $W$ is $kX$-torsion-free. Suppose that $W$ is impervious as a $kN$-module and $Se{p_{(G,Y)}}(xkG,xkY) = G$ for any element $0 \ne x \in W$ and any $G$-invariant subgroup $Y$ of $N$. Then for any finitely generated subgroup $H$ of $G$ the subgroup $N$ has an $H$-invariant polycyclic torsion-free non-abelian section. 
\end{proposition}

\begin{proof} Let $1 = {N_m} \leqslant {N_{m - 1}} \leqslant ... \leqslant {N_2} \leqslant {N_1} \leqslant {N_0} = N$ be the upper central series of the subgroup $N$. Put $K = {N_1}$ and  $M = {N_2}$ then $N' \leqslant K$. Let $L = i{s_N}(MN')$ then we have a series $M \leqslant L \leqslant K \leqslant N$ of $G$-invariant isolated subgroups of $N$ such that $N' \leqslant L$, $MN'$ is a dense subgroup of $L$ and $K/M$ is the centre of $N/M$. By Proposition 5.2, $So{c_G}N/L = N/L$ and there is a dense $G$-invariant subgroup $D \leqslant N$ such that $K \leqslant D$ and the quotient group $D/K$ is polycyclic. As $So{c_G}N/L = N/L$, we have $So{c_G}D/L = D/L$ and hence, by Lemma 2.4(iv), there is a $G$-invariant subgroup $P \leqslant D$ such that $L \leqslant P$ , the quotient group $P/L$ is polycyclic, $P \cap K = L$ and $PK$ has finite index in $D$. Then, as $D/M$ is a nilpotent torsion-free group of nilpotency class two the assertion follows from Lemma 2.1. 
\end{proof}


\section{Primitive Representations of Finitely Generated Linear Groups of Finite Rank}


\begin{theorem}
Let $G$ be an infinite finitely generated linear group of finite rank and let $k$ be a field of characteristic zero. If the group $G$ has a faithful primitive irreducible representation over the field  $k$  then $\left| \Delta (G) \right| = \infty $.
\end{theorem}
\begin{proof} If the group $G$ is polycyclic-by-finite then the result follows from \cite[Theorem A]{Harp80}. So, we may assume that the group $G$ is not polycyclic-by-finite. Suppose that the group $G$ has a faithful primitive irreducible representation over the field $k$ and let $M$ be the module of this representation. By Lemma 2.2(ii), $G$ has a finite series $L \leqslant {G_0} \leqslant G$ of normal subgroups such that $\left| G/ {G_0}\right| < \infty $, the quotient group ${G_0}/L$ is polycyclic and the subgroup $L$ is torsion-free nilpotent minimax which has no non-abelian torsion-free polycyclic $G$-invariant sections. It follows from \cite[Theorem A]{Wils70} that $M$ contains a simple $k{G_0}$-module ${M_0}$. By \cite[Theorem 4.2.5]{Tush2000}, the module ${M_0}$ is not $kL$-torsion-free and hence so is $M$. Then, as $r(G) <  \infty$, there is an isolated $G$-invariant subgroup $N \leqslant L$ such that the module $M$ is not $kN$-torsion-free but for any proper isolated $G$-invariant subgroup $X \leqslant N$ the module $M$ is $kX$-torsion-free. It follows from Lemma 2.5 that the module $M$ is $kN$-torsion. As the module $M$ is irreducible and primitive, it follows from \cite[Lemma 3.1.4(ii)]{Tush2000} that $Sep_{(G,Y)} (xkG,xkY) = G$ for any element $0 \ne x \in M$ and any $G$-invariant subgroup $Y$ of $N$. If the $kN$-module $M$ is not impervious then it follows from \cite[Proposition 3.2.3]{Tush2000} that there are an element $0 \ne x \in M$ and a proper isolated subgroup $Y \leqslant N$ such that the module $M$ is not $kY$-torsion-free and $Se{p_{(G,N)}}(xkG,xkN) \leqslant {N_G}(Y)$. Since $Sep_{(G,Y)} (xkG,xkY) = G$, we see that $G = N_G (Y)$ but it is impossible because by the choice of $N$ the module $M$ is $kY$-torsion-free for any proper isolated $G$-invariant subgroup $Y \leqslant N$. Thus, $M$ is impervious as a $kN$-module. \par  
	Suppose that $\left| {\Delta (G)} \right| < \infty $, as the subgroup $N$ is torsion-free it implies that ${\Delta(G) = 1}$. If the subgroup $N$ is abelian then by \cite[Theorem 4.2]{Tush96} there are an element $0 \ne x \in M$ and a subgroup $Y \leqslant N$ such that $xkG = (xkY){\otimes_{kY}} kG$ and ${r_0}(Y/{C_Y}(xkY)) < {r_0}(G)$. Since the module $M$ is irreducible, we have $M = xkG$ and, as the module $M$ is primitive, we see that $Y = G$. Therefore, ${r_0}(G/{C_G}(M)) < {r_0}(G)$ but it is impossible because ${C_G}(M) = 1$. 
\par  
	So, the subgroup $N$ is non-abelian. Then, by Proposition 5.3, the subgroup $N$ has a torsion-free polycyclic non-abelian $G$-invariant section but it contradicts the choice of $N$. Thus, $\left| {\Delta (G)} \right| = \infty $. 
\end{proof}




\begin{thebibliography}{30}


\bibitem{Bour} N. Bourbaki, {\small\it Elements of Mathematics: Commutative Algebra, Chapters 1-7 } 
(Springer, 1998).
\bibitem{Broo85} C. J. B. Brookes, Ideals in group rings of soluble groups of finite rank, 
{\small\it Math. Proc. Cambridge. Phil. Soc. } {\small\bf 97 } (1985) 27-49.
\bibitem{BBro85} C. J. B. Brookes and K. A. Brown, Primitive group rings and Noetherian rings of quotients, 
{\small\it Trans. Amer. Math. Soc. } {\small\bf 288 } (1985) 605-623.
\bibitem{Broo88} C. J. B. Brookes, Modules over polycyclic groups, {\small\it Proc. London Math. Soc. } 
{\small\bf 57 } (1988) 88-108. 
\bibitem{Harp77} D. L. Harper, Primitive irreducible representations of nilpotent groups, 
{\small\it Math. Proc. Camb. Phil. Soc. } {\small\bf 82 } (1977) 241-247.
\bibitem{Harp80} D. L. Harper, Primitivity in representations of polycyclic groups, 
{\small\it Math. Proc. Camb. Phil. Soc. } {\small\bf 88 } (1980) 15-31.  
\bibitem{KaMe79} M. I. Kargapolov and Yu.I. Merzljakov, {\small\it Fundamentals of the Theory of Groups},  
Graduate Texts in Math. ( Springer, 1979).
\bibitem{LeRo04} J. C. Lennox and D.J.S. Robinson, {\small\it The Theory of Infinite Soluble Groups } 
(Clarendon Press, 2004).
 \bibitem{Pass89} D. S. Passman, {\small\it Infinite Crossed Products } (Academic Press, 1989). 
\bibitem{Robi96} D. J. S. Robinson, A Course in the Theory of Groups (Springer, 1996).
 \bibitem{Szec16} F. Szechtman, Groups having a faithful irreducible representation, {\small\it J. Algebra } 
 {\small\bf 454 } (2016) 292–307.
\bibitem{SzTu17} F. Szechtman and A. Tushev, Infinite groups admitting a faithful irreducible representation, 
{\small\it J. Algebra Appl. } {\small\bf 17 } (2017) 1850005-1 - 1850005-7. 
 \bibitem{Tush90} A. V. Tushev, Irreducible representations of locally polycyclic groups over an absolute field, 
 {\small\it Ukrainian Math. } J. {\small\bf 42 } (1990) 1233-1238. 
\bibitem{Tush93} A. V. Tushev, On exact irreducible representations of locally normal groups, {\small\it Ukrainian Mat. J. } 
{\small\bf 45 } (1993) 1900-1906.
 \bibitem{Tush95} A. V. Tushev, On the primitivity of group algebras of certain classes of soluble groups of finite rank, 
 {\small\it Sb. Math. } {\small\bf 186 } (1995) 447-463.
 \bibitem{Tush96} A. V. Tushev, Spectra of conjugated ideals in group algebras of abelian groups of finite rank and control theorems, 
 {\small\it Glasgow Math. J. } {\small\bf 38 } (1996) 309-320.
\bibitem{Tush99} A. V. Tushev, Induced modules over group algebras of metabelian groups of finite rank, 
{\small\it Comm. Algebra } {\small\bf 27 } (1999) 5921-5938.  
\bibitem{Tush2000} A. V. Tushev, On the primitive representations of soluble groups of finite rank, 
 {\small\it Sb. Math. } 
 {\small\bf 191 } (2000) 117-159. 
 \bibitem{Tush02} A. V. Tushev, On Primitive Representations of Minimax Nilpotent Groups, 
 {\small\it Mathematical Notes } 
 {\small\bf 72 } (2002) 117-128.
 \bibitem{Tush12} A. V. Tushev, On the irreducible representations of soluble groups of finite rank, 
 {\small\it Asian-European J. Math. } {\small\bf 5 } (2012) 1250061-1-1250061-12.
\bibitem{Tush14} A. V. Tushev, On certain methods of studying ideals in group rings of abelian groups of finite rank, 
{\small\it Asian-European J. Math. } {\small\bf 7 } (2014) 1450065-1- 1450065-7
 \bibitem{Wils70} J. S. Wilson, Some properties of groups inherited by normal subgroups of finite index, 
 {\small\it Math. Z. } {\small\bf 144 } (1970) 19-21. 
 \bibitem{Wils88} J. S. Wilson, Soluble products of minimax groups, and nearly surjective derivations, 
 {\small\it J. Pure Appl. Algebra } {\small\bf 53 } (1988) 297-331.
 \bibitem{Wehr73} B. A. F. Wehrfritz, {\small\it Infinite linear groups } (Springer-Verlag, 1973).
 \bibitem{Wehr09} B. A. F. Wehrfritz, {\small\it Group and Ring Theoretic Properties of Polycyclic Groups}, 
 Algebra and Applications, Vol. 10 (Springer, 2009).
 \bibitem{ZaSa58} O. Zarisski and P. Samuel, {\small\it Commutative Algebra}, Vol. 1 (Van Nostrand, 1958).


\end{thebibliography}
\end{document}